\documentclass[review,12pt]{elsarticle}

\IfFileExists{mmap.sty}{\usepackage[noTeX]{mmap}}{}
\usepackage{amsfonts}
\usepackage{amsmath}
\usepackage{amssymb}
\usepackage{amsthm}
\usepackage{color}
\usepackage{setspace}
 \usepackage{mathrsfs}

\numberwithin{equation}{section}

\usepackage{subfigure}
\usepackage{graphicx}
\usepackage{paralist}
\usepackage{multirow}
\usepackage{ifthen}
\usepackage{booktabs}
\usepackage{threeparttable}
\usepackage{epsfig}
\usepackage{epstopdf}

\let\originalepsffile\epsffile
\renewcommand{\epsffile}[1]{\originalepsffile{Peimin/#1}}

\newcommand{\documentspacing}{\linespread{1.4}}
\documentspacing
\vspace{30 mm}

\begin{document}

\begin{frontmatter}

\title{ \Large Alpha-robust investment-reinsurance strategy for a mean-variance insurer under a defaultable market} 

\vskip 1.5cm

\author{  {\large Min Zhang} \\
  School of Economic Mathematics \\
Southwestern University of Finance and Economics\\ Chengdu,
611130, P.R. China \\ Email: zmin\_1022@163.com
 \vskip 0.5cm and  \vskip 0.5cm
 {\large Yong He} \\
  School of Mathematics, Physics and Data Science \\
Chongqing University of Science and Technology\\ Chongqing,
401331, P.R. China \\ Email: heyongmath@163.com
\\
\vskip 1cm Current version: Dec. 08, 2021
 \newpage
\begin{center}
{\Large Alpha-robust investment-reinsurance strategy for a mean-variance insurer under a defaultable market}
\end{center}
  }
\cortext[cor2]{Corresponding author.}

\begin{abstract}
In this paper, we consider the robust optimal reinsurance investment problem of the insurer under the $\alpha$-maxmin mean-variance criterion in the defaultable market. The financial market consists of risk-free bonds, a stock and a defaultable bond. The insurer's surplus process is described by a L\'{e}vy insurance model. From the perspective of game theory, the extended Hamilton-Jacobi-Bellman equations are established for the post-default and pre-default conditions respectively. In both cases, the closed-form expressions and corresponding value functions of the robust optimal investment reinsurance strategies are derived. Finally, numerical examples and sensitivity analysis are used to illustrate the influence of parameters on the optimal strategies.
\end{abstract}



\end{frontmatter}

\section{Introduction}
In recent years, the insurance industry has become a hot research topic in the financial field. Due to the fierce market competition, it is difficult to meet the compensation requirements of insurance companies by collecting premiums alone. In response to this problem, insurance companies generally adopt two methods. On the one hand, insurance companies make risky investments in surplus, and obtain income through investment to improve their ability to pay. On the other hand, insurance companies share part of the risk by taking the form of reinsurance. Reinsurance can effectively diversify the risk of huge claims, and investment can enable insurance companies to effectively manage surplus and maximize wealth value. Scholars have studied the optimal reinsurance investment strategies of insurance companies under different objective functions.

Some papers use the minimization of ruin probability as the objective function. For example, Hipp and Plum (2000) considered the risk process model of the compound Poisson process, and by choosing an appropriate investment strategy for the capital market index, the ruin probability of this risk process can be minimized. Young (2004) has determined the individual's best investment strategy for a specific consumption rate, and strives to minimize the possibility of lifetime bankruptcy. Promislow and Young (2005) considered investing in risky assets and purchasing proportional reinsurance to minimize the possibility of ruin. Bayraktar and Zhang (2015) determined the optimal robust investment strategy for a given consumption rate, and tried to minimize the possibility of life-time bankruptcy when she was not fully confident in the drift of risk assets. There are also most papers using utility maximization as the objective function. Yang and Zhang (2005) studied the optimal investment strategy of an insurance company with a jump-diffusion risk process. Bai and Guo (2010) considered the problem of maximizing the expected exponential utility of terminal wealth and minimizing the probability of ruin under the condition that the surplus process of the insurance company was described by Brownian motion with drift. Under the criterion of maximizing the expected utility of the index, Liang and Yuen (2016) considered the best proportional reinsurance strategy in a risk model with two dependent classes of insurance business.

Recently, more and more papers use the mean-variance criterion as the objective function. However, there is a problem of time inconsistency in the mean and variance criterion. There are two main ways to solve this problem. One is to use a pre-commitment strategy to solve this problem. For example, B{\"a}uerle (2005) studied the optimal proportion reinsurance problem of a compound Poisson risk model under the mean-variance criterion. Bai and Zhang (2008) studied the best reinsurance/new business and investment strategies for the mean-variance problem in the classic risk model and the diffusion model. Given the standard martingale approach for dealing with continuous-time portfolio selection models, Bi and Zhang (2014) considered two sub-problems to obtain the solution to the optimal problem of mean variance. Another way to solve the problem of time inconsistency is to formulate the problem in the framework of game theory. Zeng and Li (2011) studied the optimal time consistent strategy for investment reinsurance problems and investment-only problems under the mean-variance standard. Li and Li (2013) studied the optimal time consistency strategy of insurance companies based on the mean-variance criterion of state-related risk aversion. Bj{\"o}rk et al. (2014) studied the continuous-time mean-variance portfolio optimization under the framework of game theory when the dynamics of risk aversion depend on current wealth.

Although the optimal reinsurance-investment problem has been extensively studied by many scholars, only a few scholars incorporate the uncertainty of the model into the optimal reinsurance-investment problem. However, as we all know, it is difficult to accurately estimate the return on risky assets. Therefore, some scholars have advocated and studied the influence of model uncertainty on portfolio selection. Maenhout (2006) analyzed the optimal intertemporal portfolio problem of an investor who was worried about the wrong specification of the model and adhered to the robust decision-making rules when faced with the mean reversion risk premium. Liu (2010) studied the continuous-time intertemporal consumption and portfolio selection of investors with recursive preferences. Investors worry about mis-specifying models and seeking sound decision-making rules. Yi et al. (2015b) focused on the optimal portfolio selection problem with model uncertainty in financial markets involving a pair of stocks.

However, the limitation of robust utility is that it only considers the attitude of extreme vague aversion, which is not in line with reality. Therefore, some scholars consider a more general utility function: $\alpha$-maxmin expected utility. For example, Ghirardato (2004) showed how to formally model ambiguity and the reaction of decision makers in the context of general decision models. In addition, a special case of the decision model is described axiomatically, namely the "$\alpha$-maxmin" expected utility model. Inspired by the expected utility of $\alpha$-maxmin, Li et al. (2016) proposed a new mean-variance criterion, called the $\alpha$-maxmin mean-variance criterion, and applied it to the problem of reinsurance investment. On the other hand, most of the previous work only considered financial markets for two asset types: risk-free bonds and stocks, while few papers considered default securities. However, in today's financial markets, high-yield bonds are becoming increasingly attractive to investors. Although defaults rarely occur, they do exist. Therefore, the optimal portfolio of defaulted securities has become an important research topic. So far, several other studies have also tried to find the best investment strategy in the default market. Bielecki and Jang (2007) derived a closed-form solution for a representative investor who distributed his wealth among credit risk assets, non-default bank accounts and stocks. Bo et al. (2013) considered the problem of portfolio optimization in the defaultable market, and clearly derived the optimal investment and consumption strategies that maximize the indefinite expected discounted HARA utility of consumption. Zhao et al. (2016) considered the optimal investment and reinsurance of insurance companies under the mean-variance criterion in the jump diffusion risk model. The financial market includes risk-free assets, a stock and a defaultable bond. Wang et al. (2019) considered the problem of delayed reinsurance investment by insurance companies in the defaultable market under the mean-variance criterion.

In this paper, we consider an optimal investment problem, assuming that insurer can allocate his wealth in risk-free bonds, a stock and a defaultable bond. Insurer can control insurance risks by purchasing proportional reinsurance or acquiring new businesses. In addition, the surplus process of an insurance company is modeled by a spectrally negative L{\'e}vy process, while the price process of the stock follows geometric Brownian motion. Under the 
$\alpha$-maxmin mean variance criterion, we find the optimal reinsurance investment strategy with consistent time in the game theory framework. By solving the extended Hamilton-Jacobi-Bellman equations, the closed-end equilibrium reinsurance-investment strategies in the post-default and pre-default situations are derived respectively. Finally, sensitivity analysis and numerical simulation illustrate the influence of model parameters on the optimal strategies.

The rest of this paper is organized as follows: Section 2 introduces the surplus process of insurance companies and financial markets. In Section 3, we give the robust optimal reinsurance investment problem in the default market under the $\alpha$-maxmin mean-variance criterion. In Section 4, the optimal investment reinsurance strategies and the corresponding value function under the pre-default and post-default conditions are respectively obtained. Section 5 presents sensitivity analysis and numerical examples to illustrate our results. Section 6 concludes the paper.
\section{The model}
Let $\left(\Omega, \mathcal{F}, \mathbb{F}, \mathbb{P}\right)$ be a complete probability space.Let $\mathbb{F}:=\left(\mathcal{F}_{t}\right)_{t \geq 0}$ be the right-continuous, $\mathbb{P}$-complete filtration generated by two standard Brownian motions $\left\{W_{1}(t)\right\}$ and $\left\{W_{2}(t)\right\}$. We denote by $\mathbb{H}:=\left(\mathcal{H}_{t}\right)_{t \geq 0}$ the filtration of a default process $H(t)$. Let $\mathbb{G}:=\left(\mathcal{G}_{t}\right)_{t \geq 0}$ be the enlarged filtration of $\mathbb{F}$ and $\mathbb{H}$, i.e., $\mathcal{G}_{t}:=\mathcal{F}_{t} \vee \mathcal{H}_{t}$. Throughout this paper, let $\mathbb{P}$ be the real world probability measure, and $\mathbb{E}[\cdot]$ and $\operatorname{Var}[\cdot]$ be the expectation and the variance under $\mathbb{P}$, respectively. We assume that there exists a martingale probability measure $\mathbb{Q}$ (or the risk neutral measure) equivalent to measure $\mathbb{P}$. 
\subsection{Surplus process}
Suppose that the surplus process of an insurer without reinsurance and investment is simulated by a spectrally negative L{\'e}vy process defined on $\left(\Omega, \mathcal{F}, \mathbb{F}, \mathbb{P}\right)$ with dynamics
\begin{equation}
\label{201}
dR_{0}(t)=\mu_{1}dt+\sigma_{1}dW_{1}(t)-\int_{0}^{\infty} z N(dt,dz), \quad R_{0}(0)>0,
\end{equation}
where $\mu_{1}>0$ is the premium rate, $\sigma_{1}$ can be regarded as the volatility of the insurer's surplus, $\left\{W_{1}(t)\right\}_{t \geq 0}$ is a standard Brownian motion, $N(dt,dz)$ is the Poisson random measure, independent of $W_1$, and represents the number of insurance claims of size $(z,z+dz)$ within the time period $(t,t+dt)$. $\tilde{N}(dt,dz)=N(dt,dz)-v(dz)dt$ is the compensated measure of $N(dt,dz)$. And $v$ is a Lévy measure so that $\int_{0}^{\infty} z v(dz)<\infty$. We adopt the expected value principle to calculate the insurance premium. Then the premium rate satisfies $\mu_{1}=(1+\theta)\int_{0}^{\infty} z v(dz)$, where $\theta>0$ is the safety loading coefficient of the insurer.

The insurer can purchase proportional reinsurance or acquire new business to adjust the exposure to insurance risk. The proportional reinsurance or new business level level is associated with the risk exposure $\pi_{q}(t)$ at time $t$. When $\pi_{q}(t) \in[0,1]$, it corresponds to a proportional reinsurance cover. In this case, the cedent should divert part of the premium to the reinsurer at the rate of $(1+\eta)(1-\pi_{q}(t)) \int_{0}^{\infty} zv(dz)$, where $\eta$ is the safety loading of the reinsurer satisfying $\eta>\theta>0$. Then, the liability of the insurer in the period $[t,t+dt]$ reduces to $\pi_{q}(t) \int_{0}^{\infty} zN(dt,dz)dt$. In other words, the reinsurer is responsible for the remainder of the liability $(1-\pi_{q}(t)) \int_{0}^{\infty} zN(dt,dz)dt$. Note that $\pi_{q}(t) \in(1,+\infty)$ refers to acquiring new business. After taking into account such a reinsurance strategy $\pi_{q}(t)$, $t\in(0,T)$, the insurer's surplus process can be described by
\begin{equation} 
\label{202}
\begin{aligned} d R(t)&=d R_{0}(t)-(1+\eta)\left(1-\pi_{q}(t)\right) \int_{0}^{\infty} z v(d z) d t+\left(1-\pi_{q}(t)\right) \int_{0}^{\infty} z N(d t, d z) \\ &=\left[\theta-\eta+\eta \pi_{q}(t)\right] \int_{0}^{\infty} z v(d z) d t+\sigma_{1} d W_{1}(t)-\pi_{q}(t) \int_{0}^{\infty} z \tilde{N}(d t, d z). \end{aligned}
\end{equation}
\subsection{Financial market}
Consider a financial market with three available assets: a risk-free bond, a stock (or stock index) and a defaultable bond. The price process of the risk-free bond $B(t)$ under measure $\mathbb{P}$ follows
\begin{equation*}
d B(t)=r B(t) d t,
\end{equation*}
where $r>0$ denotes the risk-free interest rate. The price process of the stock under measure $\mathbb{P}$ follows
\begin{equation*}
d S(t)=S(t) \mu d t+\sigma_{2} S(t)\left(\rho d W_{1}(t)+\hat{\rho} d W_{2}(t)\right),
\end{equation*}
where $\mu \in \mathbb{R}$, $\sigma_{2}>0$, $\rho \in[-1,1]$, $\hat{\rho}=\sqrt{1-\rho^{2}}$, and $\left\{W_{2}(t)\right\}_{t \geq 0}$ is another standard Brownian motion which is independent of $W_{1}$ and $N(d t, d z)$.

Next, we describe the price dynamics of the defaultable bond under the real-world probability measure $\mathbb{P}$. Unlike the previous two securities whose dynamics have been given directly under the real world probability measure $\mathbb{P}$, the price process of the defaultable bond is first defined under the risk-neutral probability measure $\mathbb{Q}$. Then the price dynamics of the defaultable bond under $\mathbb{P}$ will be inferred from that under $\mathbb{Q}$. The bond is defaultable, we denote the default time by $\tau$, which is a nonnegative random variable in $(\Omega, \mathcal{F}, \mathbb{P})$. The default process is defined by a Poisson process $H(t):=\mathbf{1}_{\{\tau \leq t\}}$. The intensity of the jump process is denoted by $h^{\mathbb{P}}$. As shown in Bielecki and Jang (2007), the martingale jump process is given by
\begin{equation*}
M^{\mathbb{P}}(t):=H(t)-\int_{0}^{t}(1-H(u-)) h^{\mathbb{P}} d u,
\end{equation*}
which is a $(\mathbb{G}, \mathbb{P})$-martingale.

We assume that there exists  a defaultable zero-coupon bond with a maturity date $T_1$. Denote the constant loss rate by $\zeta \in[0,1]$ and the recovery rate by $1-\zeta$. The risk neutral credit spread is denoted by $\delta=h^{\mathbb{Q}} \zeta$, where $h^{\mathbb{Q}}$ is the default intensity under $\mathbb{Q}$. According to Bielecki and Jang (2007), the expected value of the defaultable bond under $\mathbb{Q}$ is given by
\begin{equation*}
p\left(t, T_{1}\right)=\mathbf{1}_{\{\tau>t\}} e^{-(r+\delta)\left(T_{1}-t\right)}+\mathbf{1}_{\{\tau \leq t\}} e^{-(r+\delta)\left(T_{1}-\tau\right)} e^{r(t-\tau)},
\end{equation*}
the dynamics of the defaultable bond price under the risk neutral
measure $\mathbb{Q}$ follows
\begin{equation*}
d p\left(t, T_{1}\right)=r p\left(t, T_{1}\right) d t-\zeta e^{-(r+\delta)\left(T_{1}-t\right)} d M^{\mathbb{Q}}(t),
\end{equation*}
where $M^{\mathbb{Q}}(t)=H(t)-\int_{0}^{t}(1-H(u-)) h^{\mathbb{Q}} d u$ is a $(\mathbb{G}, \mathbb{Q})$ martingale.

Next,we derive the price process under real world measure $\mathbb{P}$. As in Bielecki and Jang (2007), we denote by $1 / \Delta:=h^{\mathbb{Q}} / h^{\mathbb{P}}$ the constant default risk premium. Duffie and Singleton (2012) indicated that the probability of default under the risk neutral measure $\mathbb{Q}$ is higher than that under the real world probability measure $\mathbb{P}$. Therefore, we assume that $1 / \Delta \geq1$ throughout the paper.

Then by Girsanov's theorem, the price process of defaultable bond under $\mathbb{P}$ is as follows
\begin{equation}\label{203}
d p\left(t, T_{1}\right)=p\left(t-, T_{1}\right)\left[r d t+(1-H(t))(1-\Delta) \delta d t-(1-H(t-)) \zeta d M^{\mathbb{P}}(t)\right].
\end{equation}
\subsection{Wealth process}
In this paper, we assume that the insurer is allowed to purchase proportional reinsurance or acquire new business, and invest the surplus in the financial market. Let $\pi_{s}(t)$ and $\pi_{p}(t)$ be the money amount invested in the stock and the defaultable bond, respectively, and $\pi_{q}(t)$ represents the risk exposure of the insurer. The trading strategy is represented by the stochastic process $\pi=\left\{\left(\pi_{s}(u), \pi_{p}(u), \pi_{q}(u)\right)\right\}_{u \in[t, T]}$, where $T$ is the investment horizon such that $T<T_1$. Under strategy $\pi$, the insurer's wealth process $\left\{X^{\pi}(t)\right\}_{t \in[0, T]}$ follows
\begin{equation}\label{204}
\begin{aligned} d X^{\pi}(t)&=\left[r X^{\pi}(t)+(\mu-r) \pi_{s}(t)+\left(\theta-\eta+\eta \pi_{q}(t)\right) \int_{0}^{\infty} z v(d z)\right.\\ & \left.+\pi_{p}(t)(1-H(t))(1-\Delta) \delta \right] d t+\left(\sigma_{1}+\pi_{s}(t) \sigma_{2} \rho\right) d W_{1}(t)\\ &+\pi_{s}(t) \sigma_{2} \hat{\rho} d W_{2}(t)-\pi_{q}(t) \int_{0}^{\infty} z \tilde{N}(d t, d z)-\pi_{p}(t) \zeta d M^{\mathbb{P}}(t), \ X^{\pi}(0)=x_{0} \end{aligned}
\end{equation}
with $(1-H(t-)) d M^{\mathbb{P}}(t)=d M^{\mathbb{P}}(t)$ as well as the convention that $0/0=0$. This convention is need to deal with the post-default case, so that $p\left(t-, T_{1}\right)=0$ and we fix $\pi_{p}(t)=0$ afterwards.

In order to introduce the ambiguity on the insurance and financial risks, we define a set of prior probability measures as below. We call the probability distortion function $\phi:=\left(\phi_{1}(t), \phi_{2}(t), \phi_{3}(t, z)\right)_{t \in[0, T], z>0} \in \Theta$, if 
$\phi_{1}(t)$, $\phi_{2}(t)$ and $\phi_{3}(t, z)$ are deterministic functions of $t$ and $z$ and satisfy
\begin{equation*}
\begin{aligned} \exp& \left\{\int_{t}^{T} \frac{\phi_{1}(s)^{2}+\phi_{2}(s)^{2}}{2} d s\right.\\ & \left.+\int_{t}^{T} \int_{0}^{\infty}\left[\left(1-\phi_{3}(s, z)\right) \ln \left(1-\phi_{3}(s, z)\right)+\phi_{3}(s, z)\right] v(d z) d s\right\}<\infty \end{aligned}
\end{equation*}
 for any $t\in[0,T]$. Each probability distortion function $\phi \in \Theta$ is associated with a probability measure $\mathbb{Q}^{\phi} \sim \mathbb{P}$. And the Radon-Nikodym derivative process $\left.\frac{\mathrm{d} \mathbb{Q}^{\phi}}{\mathrm{d} \mathbb{P}}\right|_{\mathcal{F}(t)}:=\Lambda^{\phi}(t)$ is given by
\begin{equation*}
\begin{aligned} \Lambda^{\phi}(t) &=\exp \left\{-\int_{0}^{t} \phi_{1}(s) \mathrm{d} W_{1}(s)-\frac{1}{2} \int_{0}^{t} \phi_{1}(s)^{2} \mathrm{~d} s-\int_{0}^{t} \phi_{2}(s) \mathrm{d} W_{2}(s)-\frac{1}{2} \int_{0}^{t} \phi_{2}(s)^{2} \mathrm{~d} s\right.\\ & \left.+\int_{0}^{t} \int_{0}^{\infty} \ln \left(1-\phi_{3}(s, z)\right) \tilde{N}(\mathrm{d} s, \mathrm{d} z)\right.\\ & \left.+\int_{0}^{t} \int_{0}^{\infty}\left[\ln \left(1-\phi_{3}(s, z)\right)+\phi_{3}(s, z)\right] v(\mathrm{d} z) \mathrm{d} s\right\} . \end{aligned}
\end{equation*}
Then we define a set of prior probability measures by $\mathcal{Q}=\left\{\mathbb{Q}^{\phi}: \phi \in \Theta\right\}$. It is obvious that $\mathbb{P} \in \mathcal{Q}$ is established.

Applying Girsanov’s Theorem (see, for example, Øksendal and Sulem (2007), it is clear that
\begin{equation}\label{205}
d W_{1}^{\phi}(t)=d W_{1}(t)+\phi_{1}(t) d t\
and\
d W_{2}^{\phi}(t)=d W_{2}(t)+\phi_{2}(t) d t
\end{equation}
are $\mathbb{Q}^{\phi}$-Brownian motions, and
\begin{equation}\label{206}
\tilde{N}^{\phi}(d t, d z)=\tilde{N}(d t, d z)+\phi_{3}(t, z) v(d z) d t.
\end{equation}
Using (\ref{204})–(\ref{206}), it is easy to see that the dynamics of the surplus process $X^{\pi}(t)$ under $\mathbb{Q}^{\phi}$ is governed by
\begin{equation}\label{207}
\begin{aligned} d X^{\pi}(t)&=\left[r X^{\pi}(t)+(\mu-r)\pi_{s}(t)+\left(\theta-\eta+\eta \pi_{q}(t)\right)\int_{0}^{\infty} z v(d z)-\left(\sigma_{1}+\pi_{s}(t) \sigma_{2} \rho\right) \phi_{1}(t)\right.\\ & \left.-\pi_{s}(t) \sigma_{2} \hat{\rho} \phi_{2}(t)+\pi_{q}(t) \int_{0}^{\infty} z \phi_{3}(t, z) v(d z)+\pi_{p}(t)(1-H(t))(1-\Delta) \delta\ \right] d t\\ &+\left(\sigma_{1}+\pi_{s}(t) \sigma_{2} \rho\right) d W_{1}^{\phi}(t)+\pi_{s}(t) \sigma_{2} \hat{\rho} d W_{2}^{\phi}(t)\\ &-\pi_{q}(t) \int_{0}^{\infty} z \tilde{N}^{\phi}(d t, d z)-\pi_{p}(t) \zeta d M^{P}(t) \end{aligned}
\end{equation}
\newtheorem{defn}{Definition}[section] %
\begin{defn} A trading strategy $\pi=\left\{\left(\pi_{s}(t), \pi_{p}(t), \pi_{q}(t)\right)\right\}_{t \in[t, T]}$ is said to be admissible if it satisfies the following conditions:

(1) $\left(\left(\pi_{s}(t), \pi_{p}(t), \pi_{q}(t)\right)\right)$ is $\mathbb{Q}$-predictable;

(2) $\forall t \in[0, T], \pi_{q}(t) \geq 0$;

(3) $\mathbb{E}_{t, x}^{\underline{\phi}^{\pi}}\left[\int_{0}^{T}\left(\pi_{s}(t)^{2}+\pi_{p}(t)^{2}+\pi_{q}(t)^{2}\right) d t\right]<+\infty$ for any $(t, x) \in[0, T] \times \mathbb{R}$ and $\mathbb{E}_{t, x}^{\overline{\phi}^{\pi}}\left[\int_{0}^{T}\left(\pi_{s}(t)^{2}+\pi_{p}(t)^{2}+\pi_{q}(t)^{2}\right) d t\right]<+\infty$ for any $(t, x) \in[0, T] \times \mathbb{R}$;

(4) $\left(\pi, X^{\pi}\right)$ is the unique strong solution to the stochastic differential equation (\ref{204}).
\end{defn}
For any initial condition $(t, x, h) \in[0, T] \times \mathbb{R} \times\{0,1\}$, the corresponding set of all admissible strategies is denoted by $\Pi(t, x, h)$. Here $h$ denotes the initial state with $h=0$ and $h=1$ corresponding to the pre-default case $(\tau>t)$ and the post-default case $(\tau \leq t)$, respectively.
\section{Problem formulation}
In this section, we consider a mean-variance optimization problem for the insurer purchasing reinsurance or acquiring new business and investing in the risk-free asset, the stock and the defaultable bond. Similar to Bin et al. (2016), we define the $\alpha$-robust mean-variance criterion for a controlled surplus process $X^{\pi}$ by
\begin{equation}\label{301}
\begin{aligned} J_{\alpha}^{\pi}(t, x, h) &:=\alpha \inf _{\phi \in \Theta} \underline{J}^{\pi, \phi}(t, x, h)+\hat{\alpha} \sup _{\phi \in \Theta} \overline{J}^{\pi, \phi}(t, x, h) \\ &= \alpha \underline{J}^{\pi, \underline{\phi}^{\pi}}(t, x, h)+\hat{\alpha} \overline{J}^{\pi, \overline{\phi}^{\pi}}(t, x, h) \end{aligned}
\end{equation}
where $\alpha \in [0,1]$,
\begin{equation}\label{302}
\underline{J}^{\pi, \phi}(t, x, h)=E_{t, x, h}^{\phi}\left[X^{\pi}(T)\right]-\frac{\gamma}{2} \operatorname{Var}_{t, x, h}^{\phi}\left[X^{\pi}(T)\right]+\int_{t}^{T} h_{\beta}(\phi(s)) \mathrm{d} s,
\end{equation}
and
\begin{equation}\label{303}
\overline{J}^{\pi, \phi}(t, x, h)=E_{t, x, h}^{\phi}\left[X^{\pi}(T)\right]-\frac{\gamma}{2} \operatorname{Var}_{t, x, h}^{\phi}\left[X^{\pi}(T)\right]-\int_{t}^{T} h_{\beta}(\phi(s)) \mathrm{d} s.
\end{equation}
We write 

$\mathbb{E}_{t, x, h}^{\phi}[\cdot]=\mathbb{E}^{\phi}\left[\cdot \mid X^{\pi}(t)=x, H(t)=h\right],\operatorname{Var}_{t, x, h}^{\phi}[\cdot]=\operatorname{Var}^{\phi}\left[\cdot \mid X^{\pi}(t)=x, H(t)=h\right]$ and $\gamma>0$ is the insurer's risk aversion coefficient. $\underline{\phi}^{\pi}$ and $\overline{\phi}^{\pi}$ respectively represents
the probability distortion functions to achieve the infimum and supremum in (\ref{301}), respectively. And the penalty function is selected to
\begin{equation}\label{304}
h_{\beta}(\phi(s)):=\frac{\phi_{1}(s)^{2}}{2 \beta_{1}}+\frac{\phi_{2}(s)^{2}}{2 \beta_{2}}+\frac{\int_{0}^{\infty}\left[\left(1-\phi_{3}(s, z)\right) \ln \left(1-\phi_{3}(s, z)\right)+\phi_{3}(s, z)\right] v(\mathrm{~d} z)}{\beta_{3}}
\end{equation}
Throughout the paper, we restrict the range of $\alpha$ to 
\begin{equation}\label{305}
\frac{1}{2} \leq \alpha \leq 1
\end{equation}
In this paper, the scope of $\alpha$ is limited to $\frac{1}{2} \leq \alpha \leq 1$, $\alpha=\frac{1}{2}$ represents the ambiguity-neutral attitude, $\alpha=1$ represents the extremely ambiguity-averse attitude, and the larger the value of $\alpha$, the more ambiguity-averse attitude. In the $\alpha$-robust mean variance criterion, the deviation from the reference measure $\mathbb{P}$ is penalized by the term $\int_{t}^{T} h_{\beta}(\phi(s)) \mathrm{d} s$. The constant vector $\beta=\left(\beta_{1}, \beta_{2}, \beta_{3}\right) \in(0, \infty)^{3}$ is used to reflect the level of ambiguity about the reference measure. More specifically, $\beta_{1}, \beta_{2}$, and $\beta_{3}$ indicate the level of ambiguity about insurance premium, stock return, and insurance liability, respectively. The higher the $\beta$ value, the higher the level of ambiguity.
 A higher value of $\beta$ implies a higher level of ambiguity.
Note that
\begin{equation*}
\int_{0}^{T} h_{(1,1,1)}(\phi(s)) \mathrm{d} s=\mathbb{E}\left[\Lambda^{\phi}(T) \log \Lambda^{\phi}(T)\right]
\end{equation*}
is the so-called relative entropy of $\mathbb{Q}^{\phi}$ with respect to $\mathbb{P}$.
The main purpose of this paper is to study the $\alpha$-robust reinsurance investment problem of a class of time-consistent mean-variance insurers, that is,
\begin{equation}\label{306}
\sup _{\pi \in \Pi} J_{\alpha}^{\pi}(t, x, h),
\end{equation}
in which the equilibrium strategy are defined below.
\begin{defn} For any fixed initial state $(t, x, h) \in[0, T] \times \mathbb{R} \times\{0,1\}$, consider an admissible strategy $\pi^{*}(t, x, h)$. Choose four fixed numbers $\tilde{\pi}_{s} \in \mathbb{R}, \widetilde{\pi}_{p} \in \mathbb{R}, \widetilde{\pi}_{q} \in \mathbb{R}_{+}$ and $\epsilon \in \mathbb{R}_{+}$ and define the following straregy:
\begin{equation*}
\pi^{\epsilon}(u, \widetilde{x}, \widetilde{h}):=\left\{\begin{array}{ll}\left(\widetilde{\pi}_{s}, \widetilde{\pi}_{p}, \widetilde{\pi}_{q}\right), & \text { for }(u, \widetilde{x}, \widetilde{h}) \in[t, t+\epsilon) \times \mathbb{R} \times\{0,1\}, \\ \pi^{*}(u, \widetilde{x}, \widetilde{h}), & \text { for }(u, \widetilde{x}, \widetilde{h}) \in[t+\epsilon, T] \times \mathbb{R} \times\{0,1\}.\end{array}\right.
\end{equation*}
If
\begin{equation*}
\lim _{\epsilon \rightarrow 0} \inf \frac{J_{\alpha}^{\pi^{*}}\left(t, x, h \right)-J_{\alpha}^{\pi^{\epsilon}}\left(t, x, h\right)}{\epsilon} \geq 0
\end{equation*}
for all $\left(\widetilde{\pi}_{s}, \tilde{\pi}_{p}, \widetilde{\pi}_{q}\right) \in \mathbb{R} \times \mathbb{R} \times \mathbb{R}_{+}$ and $(t, x, h) \in[0, T] \times \mathbb{R} \times\{0,1\}$, $\pi^{*}(t, x, h)$ is called an equilibrium strategy and $J_{\alpha}^{\pi^{*}}(t, x, h)$ is the associated equilibrium value function.
\end{defn}

In this paper, we make the following integrability assumption on the L{\'e}vy measure.
\newtheorem{assumption}{Assumption}[section]
\begin{assumption} (1) If there is ambiguity to the distribution of insurance claims, that is, the function $\phi_{3}(t, z):[0, T] \times[0, \infty) \rightarrow(-\infty, 1)$ depends on $z$, we assume
\begin{equation*}
\int_{0}^{1} z v(d z)<\infty \ and \  \int_{1}^{\infty} e^{c z^{2}} v(d z)<\infty, \text{ for some } c>0.
\end{equation*}
(2) If there is no ambiguity to the distribution of insurance claims, that is, the function $\phi_{3}(t, z):[0, T] \times[0, \infty) \rightarrow(-\infty, 1)$ is independent of $z$, we assume
\begin{equation*}
v(0, \infty)<\infty \ and \ \int_{1}^{\infty} z^{2} v(\mathrm{d} z)<\infty.
\end{equation*}
\end{assumption}
The second case of Assumption 3.1 means that the jumping part of the insurance surplus process follows the compound Poisson model with jump intensity $v(0, \infty)<\infty$ under $\mathbb{P}$. As in Zeng et al (2016), this is a common simplification on the ambiguity structure for underlying jump processes.

The equilibrium strategy is time-consistent and hereafter we call the equilibrium strategy $\pi^{*}$ and the corresponding equilibrium value function $J_{\alpha}^{\pi^{*}}(t, x, h)$ the optimal time-consistent strategy and the value function for (\ref{306}), respectively. Before giving the verification theorem, we denote by
\begin{equation*}
\begin{aligned}
C^{1,2}([0, T] \times \mathbb{R}):=\{&\varphi(t, x) \mid \varphi(t, \cdot) \text{ is continuously differentiable on } [0, T],\\ &
\text{and } \varphi(\cdot, x) \text{ is twice continuously differentiable on } \mathbb{R}\},
\end{aligned}
\end{equation*}
and for any $\psi(t, x, 0), \psi(t, x, 1) \in C^{1,2}([0, T] \times \mathbb{R})$ and $\pi \in \Pi(t, x, h)$, define the infinitesimal generator
\begin{equation}\label{307}
\mathcal{A}^{\pi} \psi(t, x, h):=\left\{\begin{array}{l}\psi_{t}(t, x, 1)+[r x+\pi_{s}(\mu-r)+\left(\theta-\eta+(1+\eta) \pi_{q}\right) \int_{0}^{\infty} z v(d z)\\ \quad -\left(\sigma_{1}+\pi_{s} \sigma_{2} \rho\right) \phi_{1}(t)-\pi_{s} \sigma_{2} \hat{\rho}\phi_{2}(t)]\psi_{x}(t, x, 1)\\ \quad +\left(\frac{1}{2} \sigma_{1}^{2}+\pi_{s} \sigma_{1} \sigma_{2} \rho+\frac{1}{2} \pi_{s}^{2} \sigma_{2}^{2}\right) \psi_{x x}(t, x, 1)\\ \quad +\int_{0}^{\infty}\left(\psi\left(t, x-\pi_{q} z, 1\right)-\psi(t, x, 1)\right)\left(1-\phi_{3}(t, z)\right) v(d z),\ h=1; \\ \psi_{t}(t, x, 0)+[r x+\pi_{s}(\mu-r)+\left(\theta-\eta+(1+\eta) \pi_{q}\right) \int_{0}^{\infty} z v(d z)\\ \quad +\pi_{p} \delta-\left(\sigma_{1}+\pi_{s} \sigma_{2} \rho\right) \phi_{1}(t)-\pi_{s} \sigma_{2} \hat{\rho}\phi_{2}(t)] \psi_{x}(t, x, 0)\\ \quad +\left(\frac{1}{2} \sigma_{1}^{2}+\pi_{s} \sigma_{1} \sigma_{2} \rho+\frac{1}{2} \pi_{s}^{2} \sigma_{2}^{2}\right) \psi_{x x}(t, x, 0)\\ \quad +\int_{0}^{\infty}\left(\psi\left(t, x-\pi_{q} z, 0\right)-\psi(t, x, 0)\right)\left(1-\phi_{3}(t, z)\right) v(d z)\\ \quad +\left[\psi\left(t, x-\zeta \pi_{p}, 1\right)-\psi(t, x, 0)\right] h^{P},\ h=0.\end{array}\right.
\end{equation}
\newtheorem{thm}{\bf Theorem}[section]
\begin{thm}(Verification Theorem) For problem (\ref{306}), in the post-default case and the pre-default case, i.e., $h=1$ and $h=0$, if there exist $V(t, x, h), \underline{g}(t, x, h), \overline{g}(t, x, h) \in C^{1,2}([0, T] \times \mathbb{R})$ satisfy the following condition:

 (1) For any $(t, x, h) \in[0, T] \times \mathbb{R} \times\{0,1\}$,
\begin{equation}\label{308}
\begin{aligned}
0=\sup _{\pi \in \Pi}\left\{\alpha \inf _{\phi \in \Theta}\left[\mathcal{A}^{\pi, \phi} V(t, x, h)-\frac{\gamma}{2} \mathcal{A}^{\pi, \phi} \underline{g}^{2}(t, x, h)+\gamma g(t, x, h) \mathcal{A}^{\pi, \phi} \underline{g}(t, x, h)+h_{\beta}(\phi(t))\right]\right.\\ \left.+\hat{\alpha} \sup _{\phi \in \Theta}\left[\mathcal{A}^{\pi, \phi} V(t, x, h)-\frac{\gamma}{2} \mathcal{A}^{\pi, \phi} \overline{g}^{2}(t, x, h)+\gamma \overline{g}(t, x, h) \mathcal{A}^{\pi, \phi} \overline{g}(t, x, h)-h_{\beta}(\phi(t))\right]\right\}, \end{aligned}
\end{equation}
and $\left(\pi^{*}, \underline{\phi}^{*}, \overline{\phi}^{*}\right)$ denote the optimal values to achieve the supremum in $\pi$, infimum and supremum in $\phi$, respectively.

(2) For any $(t, x, h) \in[0, T] \times \mathbb{R} \times\{0,1\}$,
 \begin{equation}\label{309}
\left\{\begin{array}{l}V(T, x, h)=x, \\ \mathcal{A}^{\pi^{*}, \underline{\phi}^{*}} \underline{g}(t, x, h)=\mathcal{A}^{\pi^{*}, \vec{\phi}^{*}} \overline{g}(t, x, h)=0, \\ \underline{g}(T, x, h)=\overline{g}(T, x, h)=x. \end{array}\right.
\end{equation}
 
(3) For any $(t, x, h) \in[0, T] \times \mathbb{R} \times\{0,1\}$, $\pi^{*}(t)$, $\overline{\phi}^{*}(t)$, $\underline{\phi}^{*}(t)$, $\mathcal{A}^{\pi^{*}, \underline{\phi}^{*}} V(t, x, h)$, $\mathcal{A}^{\pi^{*}, \overline{\phi}^{*}} V(t, x, h)$, $\mathcal{A}^{\pi^{*}, \underline{\phi}^{*}} \underline{g}^{2}(t, x, h)$ and $\mathcal{A}^{\pi^{*}, \overline{\phi}^{*}} \overline{g}^{2}(t, x, h)$ are all deterministic functions of $t$ and independent of $x$.
 
(4) $\underline{\phi}^{*}=\underline{\phi}^{\pi^{*}}$ and $\overline{\phi}^{*}=\overline{\phi}^{\pi^{*}}$.

Then $\pi^{*}$ is the equilibrium strategy and $V(t, x, h)=J_{\alpha}^{\pi^{*}}(t, x, h)$ is the equilibrium value function to the $\alpha$-robust reinsurance-investment problem (\ref{306}). Besides, $\underline{g}(t, x, h)=\mathbb{E}_{t, x, h}^{\underline{\phi}^{*}}\left[X^{\pi^{*}}(T)\right]$ and $\overline{g}(t, x, h)=\mathbb{E}_{t, x, h}^{\overline{\phi}^{*}}\left[X^{\pi^{*}}(T)\right]$.
\end{thm}
The proof of the verification theorem can be adapted from a combination of Theorem 3.1 of Li et al [2016] and Theorem 4.1 of Bj{\"o}rk and Murgoci [2010]. We hence omit it here.

\section{Solution to the optimization problem}
In this section, we derive the optimal time-consistent reinsurance and investment strategies and the corresponding equilibrium value functions for problem (\ref{306}) in post-default case $(h=1)$ and pre-default case $(h=0)$, respectively.
\subsection{Post-default case: h=1}
In the post-default case, we have that $p\left(t, T_{1}\right)=0, \tau \leq t \leq T$. Thus $\pi_{p}(t)=0$ for $\tau \leq t \leq T$. Suppose that $V(t,x,1), \underline{g}(t,x,1), \overline{g}(t,x,1), \underline{\phi}^{*}, \overline{\phi}^{*}$ satisfying conditions (1) and (2) of Theorem 3.1. Then, with some calculations, we can rewrite (\ref{308}) as
\begin{equation}\label{401}
\begin{aligned} 0 &=\sup _{\pi \in \Pi}\left\{V_{t}+\left[r x+(\mu-r) \pi_{s}+(\theta-\eta+(1+\eta) \pi_{q}) \int_{0}^{\infty} z v(d z)\right] V_{x}\right.\\ &+\frac{1}{2}\left(\sigma_{1}^{2}+2 \pi_{s} \sigma_{1} \sigma_{2} \rho+\pi_{s}^{2} \sigma_{2}^{2}\right)\left(V_{x x}-\alpha \gamma \underline{g}_{x}^{2}-\hat{\alpha} \gamma \overline{g}_{x}^{2}\right) \\ & \left.+\alpha \inf _{\phi \in \Theta}\left\{\int_{0}^{\infty}L_{1}(V, \underline{g})\left(1-\phi_{3}(t, z)\right) v(d z) -E_{1} V_{x}+h_{\beta}(\phi)\right\}\right. \\ & \left.+\hat{\alpha} \sup _{\phi \in \Theta}\left\{\int_{0}^{\infty}L_{1}(V, \overline{g})\left(1-\phi_{3}(t, z)\right) v(d z) -E_{1} V_{x}-h_{\beta}(\phi)\right\}\right\}, \end{aligned}
\end{equation}
where
\begin{equation*}
L_{1}(V, g) =V\left(t, x-\pi_{q} z, 1\right)-V(t, x, 1)-\frac{\gamma}{2}\left(g\left(t, x-\pi_{q} z, 1\right)-g(t, x, 1)\right)^{2},
\end{equation*}
and
\begin{equation*}
E_{1}=(\sigma_{1}+\pi_{s} \sigma_{2} \rho) \phi_{1}+\pi_{s} \sigma_{2} \hat{\rho} \phi_{2}.
\end{equation*}
Applying the first-order condition on (\ref{401}) with respect to $\phi$, and obtain the following  infimum and supremum of $\phi$ respectively.
\begin{equation}\label{402}
\left\{\begin{array}{l}\underline{\phi}_{1}^{*}=\beta_{1}\left(\sigma_{1}+\pi_{s} \sigma_{2} \rho\right) V_{x}, \\ \underline{\phi}_{2}^{*}=\beta_{2} \pi_{s} \sigma_{2} \hat{\rho} V_{x}, \\ \underline{\phi}_{3}^{*}=1-\exp \left\{-\beta_{3}\left[V(t,x-\pi_{q} z, 1)-V(t, x,1)-\frac{\gamma}{2}(\underline{g}(t, x-\pi_{q} z, 1)-\underline{g}(t, x,1))^{2}\right]\right\},\end{array}\right.
\end{equation}
and
\begin{equation}\label{403}
\left\{\begin{array}{l}\overline{\phi}_{1}^{*}=-\beta_{1}\left(\sigma_{1}+\pi_{s} \sigma_{2} \rho\right) V_{x}, \\ \overline{\phi}_{2}^{*}=-\beta_{2} \pi_{s} \sigma_{2} \hat{\rho} V_{x}, \\ \overline{\phi}_{3}^{*}=1-\exp \left\{\beta_{3}\left[V(t,x-\pi_{q} z, 1)-V(t, x,1)-\frac{\gamma}{2}(\overline{g}(t, x-\pi_{q} z, 1)-\overline{g}(t, x,1))^{2}\right]\right\}.\end{array}\right.
\end{equation}
Substituting (\ref{402}) and (\ref{403}) back into (\ref{401}) yields
\begin{equation}\label{404}
\begin{aligned} 0 &=\sup _{\pi \in \Pi}\left\{V_{t}+\left[r x+(\mu-r) \pi_{s}+\left(\theta-\eta+(1+\eta) \pi_{q}\right) \int_{0}^{\infty} z v(d z)\right] V_{x}\right.\\ &+\frac{1}{2}\left(\sigma_{1}^{2}+2 \pi_{s} \sigma_{1} \sigma_{2} \rho+\pi_{s}^{2} \sigma_{2}^{2}\right)\left(V_{x x}-\alpha \gamma \underline{g}_{x}^{2}-\hat{\alpha} \gamma \overline{g}_{x}^{2}\right)\\ &+\frac{1-2 \alpha}{2}\left[\beta_{1}\left(\sigma_{1}+\pi_{s} \sigma_{2} \rho\right)^{2}+\beta_{2} \pi_{s}^{2} \sigma_{2}^{2} \hat{\rho}^{2}\right] V_{x}^{2} \\ &+\frac{\alpha}{\beta_{3}} \int_{0}^{\infty}\left[1-\exp \left\{-\beta_{3}L(V, \underline{g})\right\}\right] v(d z)\left.-\frac{\hat{\alpha}}{\beta_{3}} \int_{0}^{\infty}\left[1-\exp \left\{\beta_{3}L(V, \overline{g})\right\}\right] v(d z)\right\}. \end{aligned}
\end{equation}
By the first-order condition and differentiating with respect to $\pi_{q}(t)$ and $\pi_{s}(t)$ in (\ref{404}), respectively, we obtain $\pi_{q}^{*}(t)$ and $\pi_{s}^{*}(t)$ as follows,
\begin{equation}\label{405}
\begin{aligned} 0 &=\int_{0}^{\infty}\left\{(1+\eta) z V_{x}(t, x, 1)\right.\\ & \left.+\alpha\left[-z V_{x}\left(t, x-\pi_{q}^{*} z, 1\right)+\gamma z\left(\underline{g}\left(t, x-\pi_{q}^{*} z, 1\right)-\underline{g}(t, x, 1)\right) \underline{g}_{x}\left(t, x-\pi_{q}^{*} z, 1\right)\right]\right.\\ & \times \exp \left\{-\beta_{3}\left[V\left(t, x-\pi_{q}^{*} z, 1\right)-V(t, x, 1)-\frac{\gamma}{2}\left(\underline{g}\left(t, x-\pi_{q}^{*} z, 1\right)-\underline{g}(t, x, 1)\right)^{2}\right]\right\} \\ &+\hat{\alpha}\left[-z V_{x}\left(t, x-\pi_{q}^{*} z, 1\right)+\gamma z\left(\overline{g}\left(t, x-\pi_{q}^{*} z, 1\right)-\overline{g}(t, x, 1)\right) \overline{g}_{x}\left(t, x-\pi_{q}^{*} z, 1\right)\right] \\ & \left.\times \exp \left\{\beta_{3}\left[V\left(t, x-\pi_{q}^{*} z, 1\right)-V(t, x, 1)-\frac{\gamma}{2}\left(\overline{g}\left(t, x-\pi_{q}^{*} z, 1\right)-\overline{g}(t, x, 1)\right)^{2}\right]\right\}\right\} v(d z), \end{aligned}
\end{equation}
and
\begin{equation}\label{406}
\pi_{s}^{*}=-\frac{(\mu-r) V_{x}+\sigma_{1} \sigma_{2} \rho\left(V_{x x}-\alpha \gamma \underline{g}_{x}^{2}-\hat{\alpha} \gamma \overline{g}_{x}^{2}\right)+(1-2 \alpha) \beta_{1} \sigma_{1} \sigma_{2} \rho V_{x}^{2}}{\sigma_{2}^{2}\left[V_{x x}-\alpha \gamma \underline{g}_{x}^{2}-\hat{\alpha} \gamma \overline{g}_{x}^{2}+(1-2 \alpha)\left(\beta_{1} \rho^{2}+\beta_{2} \hat{\rho}^{2}\right) V_{x}^{2}\right]}.
\end{equation}
Inspired by the boundary conditions (\ref{309}), we assume the solutions of (\ref{404}) of the forms
\begin{equation}\label{407}
\left\{\begin{aligned} 
V(t, x, 1)&=A_{1}(t) x+B_{1}(t), \\ \underline{g}(t, x, 1)&=\underline{a}_{1}(t) x+\underline{b}_{1}(t), \\ \overline{g}(t, x, 1)&=\overline{a}_{1}(t) x+\overline{b}_{1}(t), \end{aligned}\right.
\end{equation}
where $A_{1}(t), B_{1}(t), \underline{a}_{1}(t), \underline{b}_{1}(t), \overline{a}_{1}(t), \overline{b}_{1}(t)$ are functions of $t$. By the first and the third relation of (\ref{309}), the boundary conditions are given by
\begin{equation*}
A_{1}(T)=\underline{a}_{1}(T)=\overline{a}_{1}(T)=1 \ \text{and} \  B_{1}(T)=\underline{b}_{1}(T)=\overline{b}_{1}(T)=0.
\end{equation*}
Substituting (\ref{406}) and (\ref{407}) into (\ref{404}) yields
\begin{equation}\label{408}
\begin{aligned} 0 &=A_{1}^{\prime} x+B_{1}^{\prime}+r A_{1} x+\left(\theta-\eta+(1+\eta) \pi_{q}^{*}\right) A_{1} \int_{0}^{\infty} z v(d z)-\frac{1}{2} \sigma_{1}^{2} \gamma\left(\alpha \underline{a}_{1}^{2}+\hat{\alpha} \overline{a}_{1}^{2}\right)\\ &+\frac{1-2 \alpha}{2} \beta_{1} \sigma_{1}^{2} A_{1}^{2} +\frac{\left[(\mu-r) A_{1}-\sigma_{1} \sigma_{2} \rho \gamma\left(\alpha \underline{a}_{1}^{2}+\hat{\alpha} \overline{a}_{1}^{2}\right)+(1-2 \alpha) \beta_{1} \sigma_{1} \sigma_{2} \rho A_{1}^{2}\right]^{2}}{2 \sigma_{2}^{2}\left[\alpha \gamma \underline{a}_{1}^{2}+\hat{\alpha} \gamma \overline{a}_{1}^{2}-(1-2 \alpha)\left(\beta_{1} \rho^{2}+\beta_{2} \hat{\rho}^{2}\right) A_{1}^{2}\right]} \\ &+\int_{0}^{\infty} \frac{\alpha}{\beta_{3}}\left(1-e^{\beta_{3}\left(\pi_{q}^{*} z A_{1}+\frac{\gamma}{2}\left(\pi_{q}^{*}\right)^{2} z^{2} \underline{a}_{1}^{2}\right)}\right) v(d z)\\ &-\int_{0}^{\infty} \frac{\hat{\alpha}}{\beta_{3}}\left(1-e^{-\beta_{3}\left(\pi_{q}^{*} z A_{1}+\frac{\gamma}{2}\left(\pi_{q}^{*}\right)^{2} z^{2} \overline{a}_{1}^{2}\right)}\right) v(d z). \end{aligned}
\end{equation}
Similarly, substituting (\ref{406}) and (\ref{407}) into the second relation of (\ref{309}) yields
\begin{equation}\label{409}
\begin{aligned} 0&=\underline{a}_{1}^{\prime} x+\underline{b}_{1}^{\prime}+r x \underline{a}_{1}+\left(\theta-\eta+(1+\eta) \pi_{q}^{*}\right) \underline{a}_{1} \int_{0}^{\infty} z v(d z)-\beta_{1} \sigma_{1}^{2} \underline{a}_{1} A_{1}\\ &+\left[(\mu-r) \underline{a}_{1}-2 \beta_{1} \sigma_{1} \sigma_{2} \rho \underline{a}_{1} A_{1}\right] \pi_{s}^{*}-\sigma_{2}^{2} \underline{a}_{1} A_{1}\left(\beta_{1} \rho^{2}+\beta_{2} \hat{\rho}^{2}\right)\left(\pi_{s}^{*}\right)^{2}\\ &-\pi_{q}^{*} \underline{a}_{1} \int_{0}^{\infty} z e^{\beta_{3}\left(\pi_{q}^{*} z A_{1}+\frac{\gamma}{2}\left(\pi_{q}^{*}\right)^{2} z^{2} \underline{a}_{1}^{2}\right)} v(d z), \end{aligned}
\end{equation}
and
\begin{equation}\label{410}
\begin{aligned} 0&=\overline{a}_{1}^{\prime} x+\overline{b}_{1}^{\prime}+r x \overline{a}_{1}+\left(\theta-\eta+(1+\eta) \pi_{q}^{*}\right) \overline{a}_{1} \int_{0}^{\infty} z v(d z)+\beta_{1} \sigma_{1}^{2} \overline{a}_{1} A_{1}\\ &+\left[(\mu-r) \overline{a}_{1}+2 \beta_{1} \sigma_{1} \sigma_{2} \rho \overline{a}_{1} A_{1}\right] \pi_{s}^{*}+\sigma_{2}^{2} \overline{a}_{1} A_{1}\left(\beta_{1} \rho^{2}+\beta_{2} \hat{\rho}^{2}\right)\left(\pi_{s}^{*}\right)^{2}\\ &-\pi_{q}^{*} \overline{a}_{1} \int_{0}^{\infty} z e^{-\beta_{3}\left(\pi_{q}^{*} z A_{1}+\frac{\gamma}{2}\left(\pi_{q}^{*}\right)^{2} z^{2} \overline{a}_{1}^{2}\right)} v(d z).\end{aligned}
\end{equation}
By matching the coefficients of the terms of $x$, we get
\begin{equation}\label{411}
A_{1}(t)^{\prime}+r A_{1}(t)=\underline{a}_{1}(t)^{\prime}+r \underline{a}_{1}(t)=\overline{a}_{1}(t)^{\prime}+r \overline{a}_{1}(t)=0.
\end{equation}
We use boundary conditions $A_{1}(T)=\underline{a}_{1}(T)=\overline{a}_{1}(T)=1$ to get
\begin{equation}\label{412}
A_{1}(t)=\underline{a}_{1}(t)=\overline{a}_{1}(t)=e^{r(T-t)}, \quad t \in[0, T].
\end{equation}
Substituting (\ref{412}) back into (\ref{408})-(\ref{410}), and using the boundary condition $B_{1}(T)=\underline{b}_{1}(T)=\bar{b}_{1}(T)=0$, we
can easily obtain the explicit solution of $B_{1}(t), \underline{b}_{1}(t)$, and $\overline{b}_{1}(t)$.
\begin{equation}\label{413}
\begin{aligned} B_{1}(t) &=\int_{t}^{T}\left\{\left(\theta-\eta+(1+\eta) \pi_{q}^{*}(s)\right) e^{r(T-s)} \int_{0}^{\infty} z v(d z)\right. \\ & \left.-\frac{1}{2} \gamma \sigma_{1}^{2} e^{2 r(T-s)}-\frac{2 \alpha-1}{2} \beta_{1} \sigma_{1}^{2} e^{2 r(T-s)}\right.\\ &+\frac{\left[\mu-r-\sigma_{1} \sigma_{2} \rho \gamma e^{r(T-s)}-(2 \alpha-1) \beta_{1} \sigma_{1} \sigma_{2} \rho e^{r(T-s)}\right]^{2}}{2 \sigma_{2}^{2}\left[\gamma+(2 \alpha-1)\left(\beta_{1} \rho^{2}+\beta_{2} \hat{\rho}^{2}\right)\right]} \\ &+\frac{\alpha}{\beta_{3}} \int_{0}^{\infty}\left(1-e^{\beta_{3}\left(\pi_{q}^{*}(s) z e^{r(T-s)}+\frac{\gamma}{2} \pi_{q}^{*}(s)^{2} z^{2} e^{2 r(T-s)}\right)}\right) v(d z) \\ & \left.-\frac{\hat{\alpha}}{\beta_{3}} \int_{0}^{\infty}\left(1-e^{-\beta_{3}\left(\pi_{q}^{*}(s) z e^{r(T-s)}+\frac{\gamma}{2} \pi_{q}^{*}(s)^{2} z^{2} e^{2 r(T-s)}\right)}\right) v(d z)\right\} d s, \end{aligned}
\end{equation}
and the expressions of $\underline{b}_{1}(t)$ and $\overline{b}_{1}(t)$ are given by
\begin{equation}\label{414}
\begin{aligned} \underline{b}_{1}(t) &=\int_{t}^{T}\left\{\left(\theta-\eta+(1+\eta) \pi_{q}^{*}(s)\right) e^{r(T-s)} \int_{0}^{\infty} z v(d z)-\beta_{1} \sigma_{1}^{2} e^{2 r(T-s)}\right.\\ &+\left[(\mu-r) e^{r(T-s)}-2 \beta_{1} \sigma_{1} \sigma_{2} \rho e^{2 r(T-s)}\right] \pi_{s}^{*}(s)-\sigma_{2}^{2} e^{2 r(T-s)}\left(\beta_{1} \rho^{2}+\beta_{2} \hat{\rho}^{2}\right) \pi_{s}^{*}(s)^{2} \\ & \left.-\pi_{q}^{*}(s) e^{r(T-s)} \int_{0}^{\infty} z e^{\beta_{3}\left(\pi_{q}^{*}(s) z e^{r(T-s)}+\frac{\gamma}{2} \pi_{q}^{*}(s)^{2} z^{2} e^{2 r(T-s)}\right)} v(d z)\right\} d s, \end{aligned}
\end{equation}
and
\begin{equation}\label{415}
\begin{aligned} \overline{b}_{1}(t) &=\int_{t}^{T}\left\{\left(\theta-\eta+(1+\eta) \pi_{q}^{*}(s)\right) e^{r(T-s)} \int_{0}^{\infty} z v(d z)+\beta_{1} \sigma_{1}^{2} e^{2 r(T-s)}\right.\\ &+\left[(\mu-r) e^{r(T-s)}+2 \beta_{1} \sigma_{1} \sigma_{2} \rho e^{2 r(T-s)}\right] \pi_{s}^{*}(s)+\sigma_{2}^{2} e^{2 r(T-s)}\left(\beta_{1} \rho^{2}+\beta_{2} \hat{\rho}^{2}\right) \pi_{s}^{*}(s)^{2} \\ & \left.-\pi_{q}^{*}(s) e^{r(T-s)} \int_{0}^{\infty} z e^{-\beta_{3}\left(\pi_{q}^{*}(s) z e^{r(T-s)}+\frac{\gamma}{2} \pi_{q}^{*}(s)^{2} z^{2} e^{2 r(T-s)}\right)} v(d z)\right\} d s. \end{aligned}
\end{equation}
Thus, the explicit expression of the value function is obtained. We state the equilibrium strategy for optimal control problem (\ref{306}) in the post-default case in the following theorem.
\begin{thm} For the $\alpha$-robust reinsurance-investment problem (\ref{306}), in the post-default case $(h=1)$, the optimal reinsurance and investment strategies for the insurer are as follows.

(1) The equilibrium reinsurance strategy $\pi_{q}^{*}(s)$ is determined by the equation
\begin{equation}\label{416}
\begin{aligned} 0 &=\int_{0}^{\infty}\left\{(1+\eta) z e^{r(T-t)}-\alpha\left(z e^{r(T-t)}+\gamma \pi_{q}^{*}(t) z^{2} e^{2 r(T-t)}\right) e^{\beta_{3}\left(\pi_{q}^{*}(t) z e^{(T-t)}+\frac{\gamma}{2} \pi_{q}^{*}(t)^{2} z^{2} e^{2 r(T-t)}\right)}\right.\\ & \left.-\hat{\alpha}\left(z e^{r(T-t)}+\gamma \pi_{q}^{*}(t) z^{2} e^{2 r(T-t)}\right) e^{-\beta_{3}\left(\pi_{q}^{*}(t) z e^{r(T-t)}+\frac{\gamma}{2} \pi_{q}^{*}(t)^{2} z^{2} e^{2 r(T-t)}\right)}\right\} v(d z), \end{aligned}
\end{equation}
and the equilibrium investment strategy is given by
\begin{equation}\label{417}
\pi_{s}^{*}(t)=\frac{(\mu-r) e^{-r(T-t)}-\sigma_{1} \sigma_{2} \rho \gamma-(2 \alpha-1) \beta_{1} \sigma_{1} \sigma_{2} \rho}{\sigma_{2}^{2}\left[\gamma+(2 \alpha-1)\left(\beta_{1} \rho^{2}+\beta_{2} \hat{\rho}^{2}\right)\right]}.
\end{equation}

(2) The corresponding equilibrium value function is given by
\begin{equation*}
J_{\alpha}^{\pi^{*}}(t, x, 1)=e^{r(T-t)} x+B_{1}(t),
\end{equation*}
where $B_{1}(t)$ is given by (\ref{413}).

(3) The associated probability distortion functions of the extremely ambiguity-averse measure and the extremely ambiguity-seeking measure are given respectively by
\begin{equation}\label{418}
\left\{\begin{array}{l}\underline{\phi}_{1}^{*}(t)=\beta_{1}\left(\sigma_{1}+\sigma_{2} \rho \pi_{s}^{*}(t)\right) e^{r(T-t)}, \\ \underline{\phi}_{2}^{*}(t)=\beta_{2} \sigma_{2} \hat{\rho} \pi_{s}^{*}(t) e^{r(T-t)}, \\ \underline{\phi}_{3}^{*}(t, z)=1-e^{\beta_{3}\left(\pi_{q}^{*}(t) z e^{r(T-t)}+\frac{\gamma}{2} \pi_{q}^{*}(t)^{2} z^{2} e^{2 r(T-t)}\right)},\end{array}\right.
\end{equation}
and
\begin{equation}\label{419}
\left\{\begin{array}{l}\overline{\phi}_{1}^{*}(t)=-\beta_{1}\left(\sigma_{1}+\sigma_{2} \rho \pi_{s}^{*}(t)\right) e^{r(T-t)}, \\ \overline{\phi}_{2}^{*}(t)=-\beta_{2} \sigma_{2} \hat{\rho} \pi_{s}^{*}(t) e^{r(T-t)}, \\ \overline{\phi}_{3}^{*}(t, z)=1-e^{-\beta_{3}\left(\pi_{q}^{*}(t) z e^{r(T-t)}+\frac{\gamma}{2} \pi_{q}^{*}(t)^{2} z^{2} e^{2 r(T-t)}\right)}.\end{array}\right.
\end{equation}
\end{thm}

\subsection{Pre-default case: h=0}
In the pre-default case, suppose that $V(t,x,0), \underline{g}(t,x,0), \overline{g}(t,x,0), \underline{\phi}^{*}, \overline{\phi}^{*}$ satisfying conditions (1) and (2) of Theorem 3.1. According to Theorem 3.1, (\ref{308}) can be written as
\begin{equation}\label{420}
\begin{aligned} 0 &=\sup _{\pi \in \Pi}\left\{V_{t}+\left[r x+(\mu-r) \pi_{s} +\pi_{p} \delta+(\theta-\eta+(1+\eta) \pi_{q}) \int_{0}^{\infty} z v(d z)\right] V_{x}\right.\\ &+\frac{1}{2}\left(\sigma_{1}^{2}+2 \pi_{s} \sigma_{1} \sigma_{2} \rho+\pi_{s}^{2} \sigma_{2}^{2}\right)\left(V_{x x}-\alpha \gamma \underline{g}_{x}^{2}-\hat{\alpha} \gamma \overline{g}_{x}^{2}\right)+\left[V\left(t, x-\zeta \pi_{p}, 1\right)-V(t, x, 0)\right] h^{\mathbb{P}}  \\ &-\frac{\alpha \gamma}{2}\left[\underline{g}\left(t, x-\zeta \pi_{p}, 1\right)-\underline{g}(t, x, 0)\right]^{2} h^{\mathbb{P}}-\frac{\hat{\alpha} \gamma}{2}\left[\overline{g}\left(t, x-\zeta \pi_{p}, 1\right)-\overline{g}(t, x, 0)\right]^{2} h^{\mathbb{P}} \\ & \left.+\alpha \inf _{\phi \in \Theta}\left\{\int_{0}^{\infty}L_{0}(V, \underline{g})\left(1-\phi_{3}(t, z)\right) v(d z) -E_{0} V_{x}+h_{\beta}(\phi)\right\}\right. \\ & \left.+\hat{\alpha} \sup _{\phi \in \Theta}\left\{\int_{0}^{\infty}L_{0}(V, \overline{g})\left(1-\phi_{3}(t, z)\right) v(d z) -E_{0} V_{x}-h_{\beta}(\phi)\right\}\right\}, \end{aligned}
\end{equation}
where
\begin{equation*}
L_{0}(V, g) =V\left(t, x-\pi_{q} z, 0\right)-V(t, x, 0)-\frac{\gamma}{2}\left(g\left(t, x-\pi_{q} z, 0\right)-g(t, x, 0)\right)^{2},
\end{equation*}
and
\begin{equation*}
E_{0}=(\sigma_{1}+\pi_{s} \sigma_{2} \rho) \phi_{1}+\pi_{s} \sigma_{2} \hat{\rho} \phi_{2}.
\end{equation*}
Applying the first-order condition on (\ref{420}) with respect to $\phi$, and obtain the following  infimum and supremum of $\phi$ respectively.
\begin{equation}\label{421}
\left\{\begin{array}{l}\underline{\phi}_{1}^{*}=\beta_{1}\left(\sigma_{1}+\pi_{s} \sigma_{2} \rho\right) V_{x}, \\ \underline{\phi}_{2}^{*}=\beta_{2} \pi_{s} \sigma_{2} \hat{\rho} V_{x}, \\ \underline{\phi}_{3}^{*}=1-\exp \left\{-\beta_{3}\left[V(t,x-\pi_{q} z, 0)-V(t, x,0)-\frac{\gamma}{2}(\underline{g}(t, x-\pi_{q} z, 0)-\underline{g}(t, x,0))^{2}\right]\right\},\end{array}\right.
\end{equation}
and
\begin{equation}\label{422}
\left\{\begin{array}{l}\overline{\phi}_{1}^{*}=-\beta_{1}\left(\sigma_{1}+\pi_{s} \sigma_{2} \rho\right) V_{x}, \\ \overline{\phi}_{2}^{*}=-\beta_{2} \pi_{s} \sigma_{2} \hat{\rho} V_{x}, \\ \overline{\phi}_{3}^{*}=1-\exp \left\{\beta_{3}\left[V(t,x-\pi_{q} z, 0)-V(t, x,0)-\frac{\gamma}{2}(\overline{g}(t, x-\pi_{q} z, 0)-\overline{g}(t, x,0))^{2}\right]\right\}.\end{array}\right.
\end{equation}
Substituting (\ref{421}) and (\ref{422}) back into (\ref{420}) yields
\begin{equation}\label{423}
\begin{aligned} 0 &=\sup _{\pi \in \Pi}\left\{V_{t}+\left[r x+(\mu-r) \pi_{s}+\pi_{p} \delta+\left(\theta-\eta+(1+\eta) \pi_{q}\right) \int_{0}^{\infty} z v(d z)\right] V_{x}\right.\\ &+\frac{1}{2}\left(\sigma_{1}^{2}+2 \pi_{s} \sigma_{1} \sigma_{2} \rho+\pi_{s}^{2} \sigma_{2}^{2}\right)\left(V_{x x}-\alpha \gamma \underline{g}_{x}^{2}-\hat{\alpha} \gamma \overline{g}_{x}^{2}\right)+\left[V\left(t, x-\zeta \pi_{p}, 1\right)-V(t, x, 0)\right] h^{\mathbb{P}} \\ &-\frac{\alpha \gamma}{2}\left[\underline{g}\left(t, x-\zeta \pi_{p}, 1\right)-\underline{g}(t, x, 0)\right]^{2} h^{\mathbb{P}}-\frac{\hat{\alpha} \gamma}{2}\left[\overline{g}\left(t, x-\zeta \pi_{p}, 1\right)-\overline{g}(t, x, 0)\right]^{2} h^{\mathbb{P}} \\ &+\frac{1-2 \alpha}{2}\left[\beta_{1}\left(\sigma_{1}+\pi_{s} \sigma_{2} \rho\right)^{2}+\beta_{2} \pi_{s}^{2} \sigma_{2}^{2} \hat{\rho}^{2}\right] V_{x}^{2} \\ &+\frac{\alpha}{\beta_{3}} \int_{0}^{\infty}\left[1-\exp \left\{-\beta_{3}L_{0}(V, \underline{g})\right\}\right] v(d z)\left.-\frac{\hat{\alpha}}{\beta_{3}} \int_{0}^{\infty}\left[1-\exp \left\{\beta_{3}L_{0}(V, \overline{g})\right\}\right] v(d z)\right\}. \end{aligned}
\end{equation}
By differentiating with respect to $\pi_{s}(t)$, $\pi_{p}(t)$ and $\pi_{q}(t)$, we have
\begin{equation}\label{424}
\pi_{s}^{*}=-\frac{(\mu-r) V_{x}+\sigma_{1} \sigma_{2} \rho\left(V_{x x}-\alpha \gamma \underline{g}_{x}^{2}-\hat{\alpha} \gamma \overline{g}_{x}^{2}\right)+(1-2 \alpha) \beta_{1} \sigma_{1} \sigma_{2} \rho V_{x}^{2}}{\sigma_{2}^{2}\left[V_{x x}-\alpha \gamma \underline{g}_{x}^{2}-\hat{\alpha} \gamma \overline{g}_{x}^{2}+(1-2 \alpha)\left(\beta_{1} \rho^{2}+\beta_{2} \hat{\rho}^{2}\right) V_{x}^{2}\right]},
\end{equation}
\begin{equation}\label{425}
\begin{aligned} 0 &=\delta V_{x}-\zeta h^{\mathbb{P}} V_{x}\left(t, x-\zeta \pi_{p}^{*}, 1\right)+\alpha \gamma \zeta h^{\mathbb{P}} \underline{g}_{x}\left(t, x-\zeta \pi_{p}^{*}, 1\right)\left[\underline{g}\left(t, x-\zeta \pi_{p}^{*}, 1\right)-\underline{g}(t, x, 0)\right] \\ &+\hat{\alpha} \gamma \zeta h^{\mathbb{P}} \overline{g}_{x}\left(t, x-\zeta \pi_{p}^{*}, 1\right)\left[\overline{g}\left(t, x-\zeta \pi_{p}^{*}, 1\right)-\overline{g}(t, x, 0)\right], \end{aligned}
\end{equation}
\begin{equation}\label{426}
\begin{aligned} 0 &=\int_{0}^{\infty}\left\{(1+\eta) z V_{x}(t, x, 0)\right.\\ & \left.+\alpha\left[-z V_{x}\left(t, x-\pi_{q}^{*} z, 0\right)+\gamma z\left(\underline{g}\left(t, x-\pi_{q}^{*} z, 0\right)-\underline{g}(t, x, 0)\right) \underline{g}_{x}\left(t, x-\pi_{q}^{*} z, 0\right)\right]\right.\\ & \times \exp \left\{-\beta_{3}\left[V\left(t, x-\pi_{q}^{*} z, 0\right)-V(t, x, 0)-\frac{\gamma}{2}\left(\underline{g}\left(t, x-\pi_{q}^{*} z, 0\right)-\underline{g}(t, x, 0)\right)^{2}\right]\right\} \\ &+\hat{\alpha}\left[-z V_{x}\left(t, x-\pi_{q}^{*} z, 0\right)+\gamma z\left(\overline{g}\left(t, x-\pi_{q}^{*} z, 0\right)-\overline{g}(t, x, 0)\right) \overline{g}_{x}\left(t, x-\pi_{q}^{*} z, 0\right)\right] \\ & \left.\times \exp \left\{\beta_{3}\left[V\left(t, x-\pi_{q}^{*} z, 0\right)-V(t, x, 0)-\frac{\gamma}{2}\left(\overline{g}\left(t, x-\pi_{q}^{*} z, 0\right)-\overline{g}(t, x, 0)\right)^{2}\right]\right\}\right\} v(d z). \end{aligned}
\end{equation}
Similarly, we conjecture $V(t, x, 0)$, $\underline{g}(t, x, 0)$ and $\overline{g}(t, x, 0)$ are the following parametric forms:
\begin{equation}\label{427}
\left\{\begin{aligned} 
V(t, x, 0)&=A_{0}(t) x+B_{0}(t), \\ \underline{g}(t, x, 0)&=\underline{a}_{0}(t) x+\underline{b}_{0}(t), \\ \overline{g}(t, x, 0)&=\overline{a}_{0}(t) x+\overline{b}_{0}(t), \end{aligned}\right.
\end{equation}
where $A_{0}(t), B_{0}(t), \underline{a}_{0}(t), \underline{b}_{0}(t), \overline{a}_{0}(t), \overline{b}_{0}(t)$ are functions of $t$. By the first and the third relation of (\ref{309}), the boundary conditions are given by
\begin{equation*}
A_{0}(T)=\underline{a}_{0}(T)=\overline{a}_{0}(T)=1 \ \text{and} \  B_{0}(T)=\underline{b}_{0}(T)=\overline{b}_{0}(T)=0.
\end{equation*}
Substituting (\ref{424}) and (\ref{427}) into (\ref{423}) yields
\begin{equation}\label{428}
\begin{aligned} 0 &=A_{0}^{\prime} x+B_{0}^{\prime}+r A_{0} x+\left(\theta-\eta+(1+\eta) \pi_{q}^{*}\right) A_{0} \int_{0}^{\infty} z v(d z)-\frac{1}{2} \sigma_{1}^{2} \gamma\left(\alpha \underline{a}_{0}^{2}+\hat{\alpha} \overline{a}_{0}^{2}\right) \\ &+\frac{1-2 \alpha}{2} \beta_{1} \sigma_{1}^{2} A_{0}^{2}+\pi_{p}^{*} \delta A_{0}+\frac{\left[(\mu-r) A_{0}-\sigma_{1} \sigma_{2} \rho \gamma\left(\alpha \underline{a}_{0}^{2}+\hat{\alpha} \overline{a}_{0}^{2}\right)+(1-2 \alpha) \beta_{1} \sigma_{1} \sigma_{2} \rho A_{0}^{2}\right]^{2}}{2 \sigma_{2}^{2}\left[\alpha \gamma \underline{a}_{0}^{2}+\hat{\alpha} \gamma \overline{a}_{0}^{2}-(1-2 \alpha)\left(\beta_{1} \rho^{2}+\beta_{2} \hat{\rho}^{2}\right) A_{0}^{2}\right]} \\ &+h^{\mathbb{P}}\left[A_{1}\left(x-\zeta \pi_{p}^{*}\right)-A_{0} x+B_{1}-B_{0}\right]-\frac{\alpha \gamma}{2}\left[\underline{a}_{1}\left(x-\zeta \pi_{p}^{*}\right)-\underline{a}_{0} x+\underline{b}_{1}-\underline{b}_{0}\right]^{2} h^{\mathbb{P}} \\ &-\frac{\hat{\alpha} \gamma}{2}\left[\overline{a}_{1}\left(x-\zeta \pi_{p}^{*}\right)-\overline{a}_{0} x+\overline{b}_{1}-\overline{b}_{0}\right]^{2} h^{\mathbb{P}} \\ &+\int_{0}^{\infty} \frac{\alpha}{\beta_{3}}\left(1-e^{\beta_{3}\left(\pi_{q}^{*}z A_{0}+\frac{\gamma}{2}\left(\pi_{q}^{*}\right)^{2}z^{2} \underline{a}_{0}^{2}\right)}\right)-\frac{\hat{\alpha}}{\beta_{3}}\left(1-e^{-\beta_{3}\left(\pi_{q}^{*}z A_{0}+\frac{\gamma}{2}\left(\pi_{q}^{*}\right)^{2} z^{2} \overline{a}_{0}^{2}\right)}\right) v(d z). \end{aligned}
\end{equation}
Similarly, substituting (\ref{424}) and (\ref{427}) into the second relation of (\ref{309}) yields
\begin{equation}\label{429}
\begin{aligned} 0 &=\underline{a}_{0}^{\prime} x+\underline{b}_{0}^{\prime}+r x \underline{a}_{0}+\left(\theta-\eta+(1+\eta) \pi_{q}^{*}\right) \underline{a}_{0} \int_{0}^{\infty} z v(d z)-\beta_{1} \sigma_{1}^{2} \underline{a}_{0} A_{0}\\ &+\left[(\mu-r) \underline{a}_{0}-2 \beta_{1} \sigma_{1} \sigma_{2} \rho \underline{a}_{0} A_{0}\right] \pi_{s}^{*} -\sigma_{2}^{2} \underline{a}_{0} A_{0}\left(\beta_{1} \rho^{2}+\beta_{2} \hat{\rho}^{2}\right)\left(\pi_{s}^{*}\right)^{2} \\ &-\pi_{q}^{*} \underline{a}_{0} \int_{0}^{\infty} z e^{\beta_{3}\left(\pi_{q}^{*}z A_{0}+\frac{\gamma}{2}\left(\pi_{q}^{*}\right)^{2} z^{2} \underline{a}_{0}^{2}\right)} v(d z) \\ &+\pi_{p}^{*} \delta \underline{a}_{0}+\left[\underline{a}_{1}\left(x-\zeta \pi_{p}^{*}\right)-\underline{a}_{0} x+\underline{b}_{1}-\underline{b}_{0}\right] h^{\mathbb{P}}, \end{aligned}
\end{equation}
and
\begin{equation}\label{430}
\begin{aligned} 0&=\overline{a}_{0}^{\prime} x+\overline{b}_{0}^{\prime}+r x \overline{a}_{0}+\left(\theta-\eta+(1+\eta) \pi_{q}^{*}\right) \overline{a}_{0} \int_{0}^{\infty} z v(d z)+\beta_{1} \sigma_{1}^{2} \overline{a}_{0} A_{0}\\ &+\left[(\mu-r) \overline{a}_{0}+2 \beta_{1} \sigma_{1} \sigma_{2} \rho \overline{a}_{0} A_{0}\right] \pi_{s}^{*} +\sigma_{2}^{2} \overline{a}_{0} A_{0}\left(\beta_{1} \rho^{2}+\beta_{2} \hat{\rho}^{2}\right)\left(\pi_{s}^{*}\right)^{2} \\ &-\pi_{q}^{*} \overline{a}_{0} \int_{0}^{\infty} z e^{-\beta_{3}\left(\pi_{q}^{*}z A_{0}+\frac{\gamma}{2}\left(\pi_{q}^{*}\right)^{2} z^{2} \overline{a}_{0}^{2}\right)} v(d z)\\ &+\pi_{p}^{*} \delta \overline{a}_{0}+\left[\overline{a}_{1}\left(x-\zeta \pi_{p}^{*}\right)-\overline{a}_{0} x+\overline{b}_{1}-\overline{b}_{0}\right] h^{\mathbb{P}}. \end{aligned}
\end{equation}
Matching the coefficients of the terms of $x$ and using boundary conditions $A_{0}(T)=\underline{a}_{0}(T)=\overline{a}_{0}(T)=1$ leads to
\begin{equation}\label{431}
A_{0}(t)=\underline{a}_{0}(t)=\overline{a}_{0}(t)=e^{r(T-t)}, \quad t \in[0, T].
\end{equation}
Substituting (\ref{431}) back into (\ref{428})-(\ref{430}), and using the boundary condition $B_{0}(T)=\underline{b}_{0}(T)=\bar{b}_{0}(T)=0$, we can easily obtain the explicit solution of $B_{0}(t), \underline{b}_{0}(t)$, and $\overline{b}_{0}(t)$.
\begin{equation}\label{432}
\begin{aligned} B_{0}(t) &=e^{-h^{\mathbb{P}}(T-t)} \int_{t}^{T}e^{h^{\mathbb{P}}(T-s)} \left\{\left(\theta-\eta+(1+\eta) \pi_{q}^{*}(s)\right) e^{r(T-s)} \int_{0}^{\infty} z v(d z)\right. \\ & \left.-\frac{1}{2} \gamma \sigma_{1}^{2} e^{2 r(T-s)}-\frac{2 \alpha-1}{2} \beta_{1} \sigma_{1}^{2} e^{2 r(T-s)}+\pi_{p}^{*}(s) \delta e^{r(T-s)}+h^{\mathbb{P}}\left[-e^{r(T-s)} \pi_{p}^{*}(s) \zeta+B_{1}\right]\right. \\ &-\frac{\alpha \gamma}{2}\left[-e^{r(T-s)} \pi_{p}^{*}(s) \zeta+\underline{b}_{1}-\underline{b}_{0}\right]^{2} h^{\mathbb{P}}-\frac{\hat{\alpha} \gamma}{2}\left[-e^{r(T-s)} \pi_{p}^{*}(s) \zeta+\overline{b}_{1}-\overline{b}_{0}\right]^{2} h^{\mathbb{P}} \\ &+\frac{\left[\mu-r-\sigma_{1} \sigma_{2} \rho \gamma e^{r(T-s)}-(2 \alpha-1) \beta_{1} \sigma_{1} \sigma_{2} \rho e^{r(T-s)}\right]^{2}}{2 \sigma_{2}^{2}\left[\gamma+(2 \alpha-1)\left(\beta_{1} \rho^{2}+\beta_{2} \hat{\rho}^{2}\right)\right]}\\ &+\frac{\alpha}{\beta_{3}} \int_{0}^{\infty}\left(1-e^{\beta_{3}\left(\pi_{q}^{*}(s) z e^{r(T-s)}+\frac{\gamma}{2} \pi_{q}^{*}(s)^{2} z^{2} e^{2 r(T-s)}\right)}\right) v(d z) \\ & \left.-\frac{\hat{\alpha}}{\beta_{3}} \int_{0}^{\infty}\left(1-e^{-\beta_{3}\left(\pi_{q}^{*}(s) z e^{r(T-s)}+\frac{\gamma}{2} \pi_{q}^{*}(s)^{2} z^{2} e^{2 r(T-s)}\right)}\right) v(d z)\right\} d s, \end{aligned}
\end{equation}
and the expressions of $\underline{b}_{0}(t)$ and $\overline{b}_{0}(t)$ are given by
\begin{equation}\label{433}
\begin{aligned} \underline{b}_{0}(t) &=e^{-h^{\mathbb{P}}(T-t)} \int_{t}^{T}e^{h^{\mathbb{P}}(T-s)}\left\{\left(\theta-\eta+(1+\eta) \pi_{q}^{*}(s)\right) e^{r(T-s)} \int_{0}^{\infty} z v(d z)-\beta_{1} \sigma_{1}^{2} e^{2 r(T-s)}\right.\\ &+\left[(\mu-r) e^{r(T-s)}-2 \beta_{1} \sigma_{1} \sigma_{2} \rho e^{2 r(T-s)}\right] \pi_{s}^{*}(s)-\sigma_{2}^{2} e^{2 r(T-s)}\left(\beta_{1} \rho^{2}+\beta_{2} \hat{\rho}^{2}\right) \pi_{s}^{*}(s)^{2} \\ & \left.-\pi_{q}^{*}(s) e^{r(T-s)} \int_{0}^{\infty} z e^{\beta_{3}\left(\pi_{q}^{*}(s) z e^{r(T-s)}+\frac{\gamma}{2} \pi_{q}^{*}(s)^{2} z^{2} e^{2 r(T-s)}\right)} v(d z) \right.\\ & \left.+\left[-e^{r(T-s)} \pi_{p}^{*}(s) \zeta+\underline{b}_{1}\right] h^{\mathbb{P}}+\pi_{p}^{*}(s) \delta e^{r(T-s)}\right\} d s, \end{aligned}
\end{equation}
and
\begin{equation}\label{434}
\begin{aligned} \overline{b}_{0}(t) &=e^{-h^{\mathbb{P}}(T-t)} \int_{t}^{T}e^{h^{\mathbb{P}}(T-s)}\left\{\left(\theta-\eta+(1+\eta) \pi_{q}^{*}(s)\right) e^{r(T-s)} \int_{0}^{\infty} z v(d z)+\beta_{1} \sigma_{1}^{2} e^{2 r(T-s)}\right.\\ &+\left[(\mu-r) e^{r(T-s)}+2 \beta_{1} \sigma_{1} \sigma_{2} \rho e^{2 r(T-s)}\right] \pi_{s}^{*}(s)+\sigma_{2}^{2} e^{2 r(T-s)}\left(\beta_{1} \rho^{2}+\beta_{2} \hat{\rho}^{2}\right) \pi_{s}^{*}(s)^{2} \\ & \left.-\pi_{q}^{*}(s) e^{r(T-s)} \int_{0}^{\infty} z e^{-\beta_{3}\left(\pi_{q}^{*}(s) z e^{r(T-s)}+\frac{\gamma}{2} \pi_{q}^{*}(s)^{2} z^{2} e^{2 r(T-s)}\right)} v(d z) \right.\\ & \left.+\left[-e^{r(T-s)} \pi_{p}^{*}(s) \zeta+\overline{b}_{1}\right] h^{\mathbb{P}}+\pi_{p}^{*}(s) \delta e^{r(T-s)}\right\} d s. \end{aligned}
\end{equation}

The above results are summarized in the following theorem.
\begin{thm} For the $\alpha$-robust reinsurance-investment problem (\ref{306}), in the pre-default case $(h=0)$, the optimal reinsurance and investment strategies for the insurer are as follows.

(1) The equilibrium reinsurance strategy $\pi_{q}^{*}(t)$ is determined by the equation
\begin{equation}\label{435}
\begin{aligned} 0 &=\int_{0}^{\infty}\left\{(1+\eta) z e^{r(T-t)}-\alpha\left(z e^{r(T-t)}+\gamma \pi_{q}^{*}(t) z^{2} e^{2 r(T-t)}\right) e^{\beta_{3}\left(\pi_{q}^{*}(t) z e^{(T-t)}+\frac{\gamma}{2} \pi_{q}^{*}(t)^{2} z^{2} e^{2 r(T-t)}\right)}\right.\\ & \left.-\hat{\alpha}\left(z e^{r(T-t)}+\gamma \pi_{q}^{*}(t) z^{2} e^{2 r(T-t)}\right) e^{-\beta_{3}\left(\pi_{q}^{*}(t) z e^{r(T-t)}+\frac{\gamma}{2} \pi_{q}^{*}(t)^{2} z^{2} e^{2 r(T-t)}\right)}\right\} v(d z), \end{aligned}
\end{equation}
and the equilibrium investment strategy is given by
\begin{equation}\label{436}
\pi_{s}^{*}(t)=\frac{(\mu-r) e^{-r(T-t)}-\sigma_{1} \sigma_{2} \rho \gamma-(2 \alpha-1) \beta_{1} \sigma_{1} \sigma_{2} \rho}{\sigma_{2}^{2}\left[\gamma+(2 \alpha-1)\left(\beta_{1} \rho^{2}+\beta_{2} \hat{\rho}^{2}\right)\right]},
\end{equation}
and
\begin{equation}\label{437}
\begin{aligned}
0&=\delta e^{r(T-t)}-\zeta h^{\mathbb{P}} e^{r(T-t)}+\alpha \gamma \zeta h^{\mathbb{P}} e^{r(T-t)}\left(-\zeta \pi_{p}^{*}(t) e^{r(T-t)}+\underline{b}_{1}(t)-\underline{b}_{0}(t)\right)\\ &+\hat{\alpha} \gamma \zeta h^{\mathbb{P}} e^{r(T-t)}\left(-\zeta \pi_{p}^{*}(t) e^{r(T-t)}+\overline{b}_{1}(t)-\overline{b}_{0}(t)\right),
\end{aligned}
\end{equation}
where $\underline{b}_{1}(t)$, $\overline{b}_{1}(t)$, $\underline{b}_{0}(t)$ and $\underline{b}_{0}(t)$ are given by (\ref{414}), (\ref{415}), (\ref{433}) and (\ref{434}), respectively.

(2) The corresponding equilibrium value function is given by
\begin{equation*}
J_{\alpha}^{\pi^{*}}(t, x, 0)=e^{r(T-t)} x+B_{0}(t),
\end{equation*}
where $B_{0}(t)$ is given by (\ref{432}).

(3) The associated probability distortion functions of the extremely ambiguity-averse measure and the extremely ambiguity-seeking measure are given respectively by
\begin{equation}\label{438}
\left\{\begin{array}{l}\underline{\phi}_{1}^{*}(t)=\beta_{1}\left(\sigma_{1}+\sigma_{2} \rho \pi_{s}^{*}(t)\right) e^{r(T-t)}, \\ \underline{\phi}_{2}^{*}(t)=\beta_{2} \sigma_{2} \hat{\rho} \pi_{s}^{*}(t) e^{r(T-t)}, \\ \underline{\phi}_{3}^{*}(t, z)=1-e^{\beta_{3}\left(\pi_{q}^{*}(t) z e^{r(T-t)}+\frac{\gamma}{2} \pi_{q}^{*}(t)^{2} z^{2} e^{2 r(T-t)}\right)},\end{array}\right.
\end{equation}
and
\begin{equation}\label{439}
\left\{\begin{array}{l}\overline{\phi}_{1}^{*}(t)=-\beta_{1}\left(\sigma_{1}+\sigma_{2} \rho \pi_{s}^{*}(t)\right) e^{r(T-t)}, \\ \overline{\phi}_{2}^{*}(t)=-\beta_{2} \sigma_{2} \hat{\rho} \pi_{s}^{*}(t) e^{r(T-t)}, \\ \overline{\phi}_{3}^{*}(t, z)=1-e^{-\beta_{3}\left(\pi_{q}^{*}(t) z e^{r(T-t)}+\frac{\gamma}{2} \pi_{q}^{*}(t)^{2} z^{2} e^{2 r(T-t)}\right)}.\end{array}\right.
\end{equation}
\end{thm}

\section{Numerical analysis}
In this section, we present some numerical examples to illustrate the effects of model parameters on the results derived in the previous section. For numerical illustrations,unless otherwise stated, the basic parameters are follows: $\alpha=0.8$, $\gamma=0.5$, $\theta=0.1$, $\eta=0.2$, $\beta_1=1$, $\beta_2=3$, $\beta_3=0.1$, $\mu=0.1$, $r=0.05$, $\lambda=1$, $\sigma_1=0.5$, $\sigma_2=0.2$, $\rho=-0.5$, $\delta=0.01$, $\zeta=0.5$, $h^{\mathbb{P}}=0.002$, $t=0$ and $T=10$.

\subsection{Analysis of the optimal reinsurance strategy}
Assuming that the total insurance claims follow a compound Poisson structure, that is, the dynamics of the insurance surplus process without reinsurance and investment given in (\ref{201}) can be rewritten as
\begin{equation*}
R_{0}(t)=R_{0}(0)+\mu_1 t+\sigma_{1} W_{1}(t)-\sum_{i=1}^{N(t)} Z_{i},
\end{equation*}
where $\{N(t)\}_{t \geq 0}$ is a homogeneous Poisson process with intensity $\lambda>0$, and $\left\{Z_{i}\right\}_{i \in \mathbb{N}}$ is a positive, independent and identically distributed random variable sequence. In addition, ${Z_{i}}$ is assumed to follow the truncated normal distribution on $(0, \infty)$, with the parameters $\mu_{Z}=1$ and $\sigma_{Z}=0.1$. Therefore, the corresponding L{\'e}vy measure is given by
\begin{equation*}
v(d z)=\lambda \frac{\frac{1}{\sigma_{Z} \sqrt{2 \pi}} \exp \left(-\frac{\left(z-\mu_{Z}\right)^{2}}{2 \sigma_{Z}^{2}}\right)}{1-\Phi\left(\frac{-\mu_{Z}}{\sigma_{Z}}\right)} d z, \quad z>0,
\end{equation*}
which clearly satisfies part (1) of Assumption 3.1. Since the equation (\ref{416}) and the equation (\ref{435}) are the same, we bring other parameters into the equation (\ref{435}) to solve.

\begin{figure}[htpp]
    \centering
    \subfigure[Effect of $\alpha$ on $\pi_q^*$]{
    \includegraphics[width=7cm]{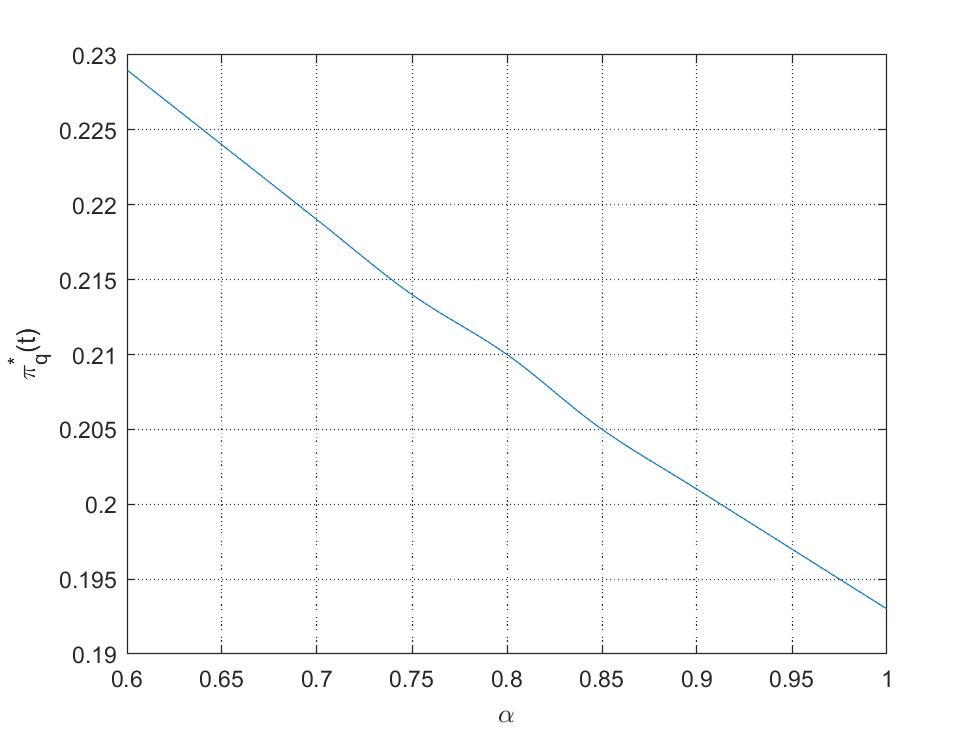}
    }
    \quad
    \subfigure[Effect of $\beta_3$ on $\pi_q^*$]{
    \includegraphics[width=7cm]{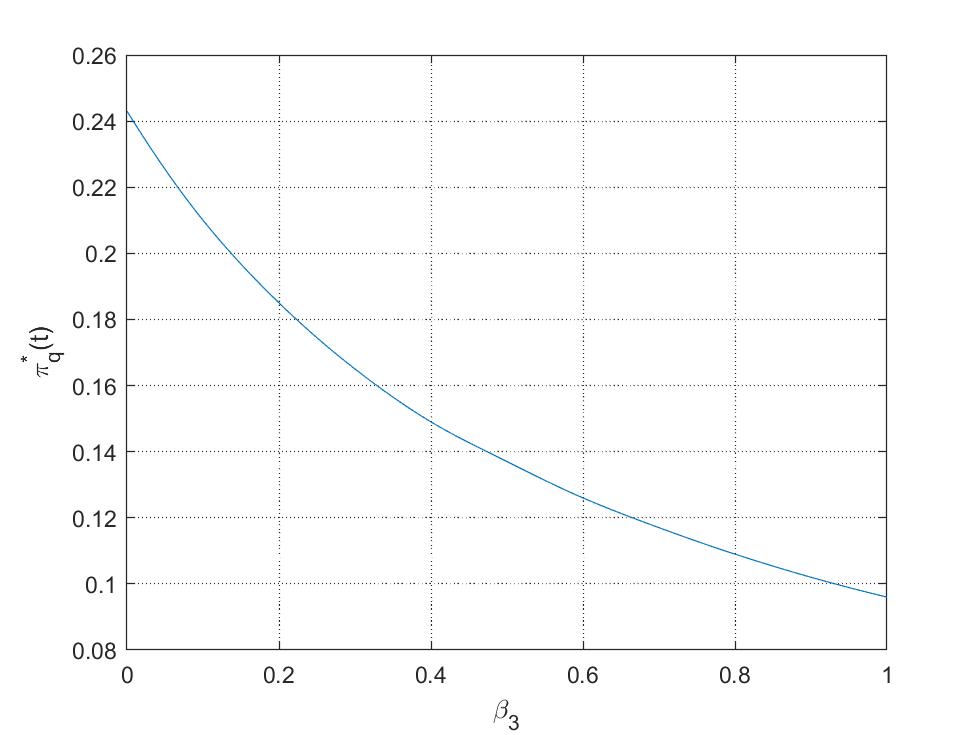}
    }
    \caption{Effects of parameters $\alpha$ and $\beta_3$ on $\pi_q^*$}
    \label{fig:my_label_1}
\end{figure}

\begin{figure}[htpp]
    \centering
    \subfigure[Effect of $\gamma$ on $\pi_q^*$]{
    \includegraphics[width=7cm]{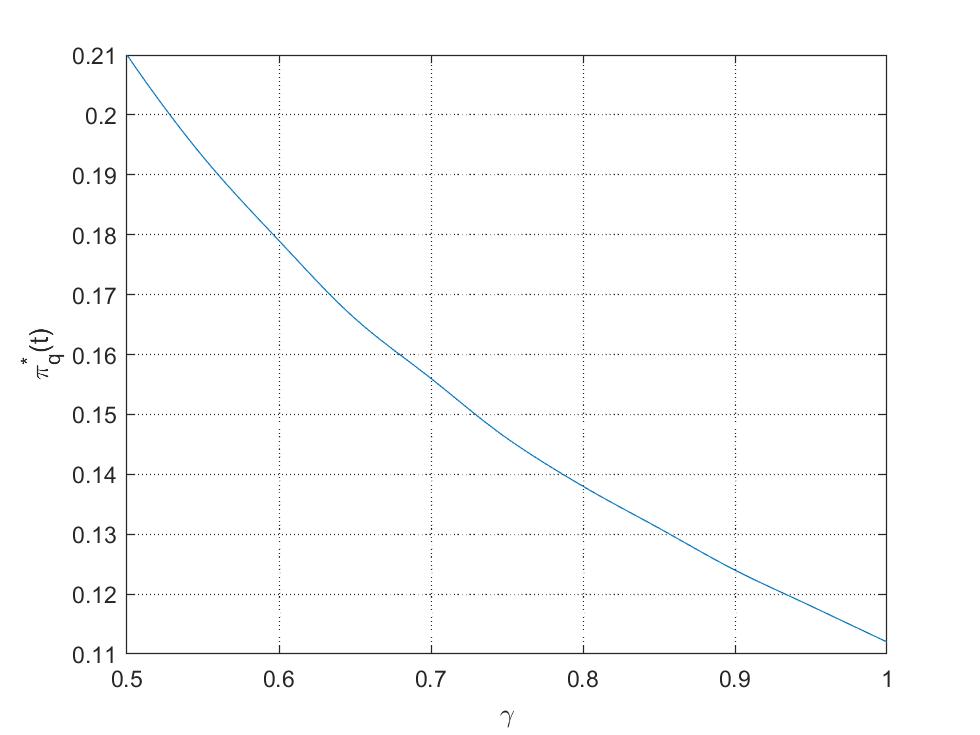}
    }
    \quad
    \subfigure[Effect of $t$ on $\pi_q^*$]{
    \includegraphics[width=7cm]{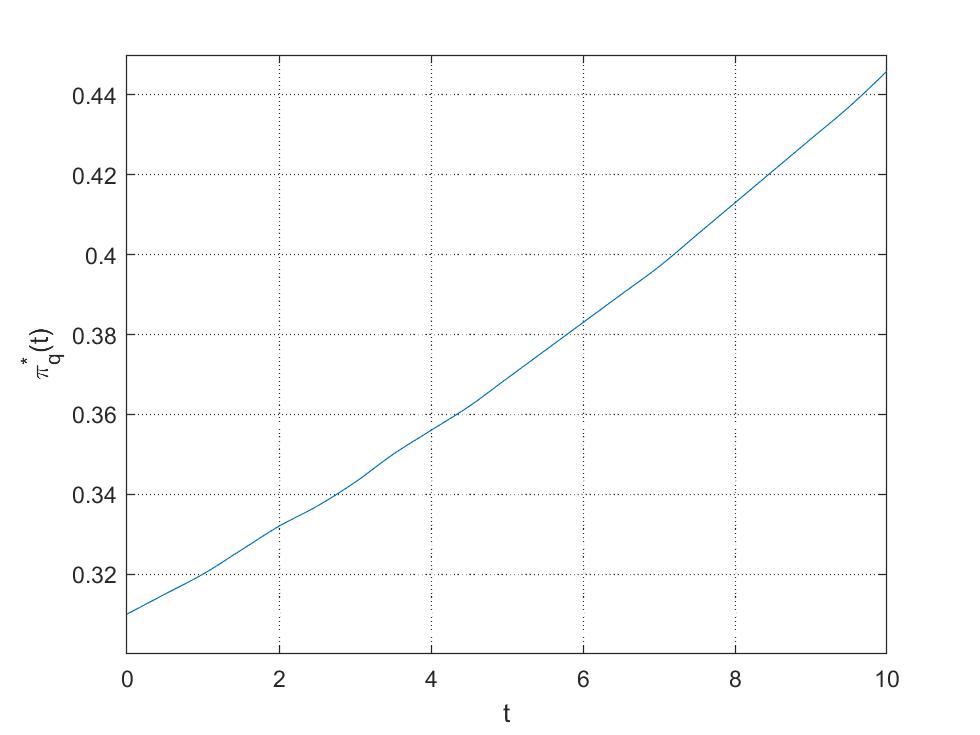}
    }
    \caption{Effects of parameters $\gamma$ and $t$ on $\pi_q^*$}
    \label{fig:my_label_2}
\end{figure}
Figure 1 and Figure 2 show that the optimal reinsurance strategy $\pi_q^*(t)$ decreases with ambiguity attitude $\alpha$, the level of ambiguity towards insurance liability $\beta_3$ and risk coefficient $\gamma$, while increases with the current time $t$. As $\alpha$ and $\gamma$ increase, the insurer is more ambiguity averse and risk averse, the insurer will be more conservative to undertake risk and thus purchase more reinsurance or acquire less new business. As $\beta_3$ increase, the insurance liability is more ambiguous, then to manage the risk exposure, the insurer will purchase more reinsurance or acquire less new business. In addition, as t increases, the insurer will preserve more insurance business by purchasing less reinsurance or acquiring more new business.

\subsection{Analysis of the optimal investment strategy}
In this subsection, we provide some sensitivity analyses on the effect of the parameters on the optimal investment strategies $\pi_s^*(t)$ and $\pi_p^*(t)$.

\begin{figure}[htpp]
    \centering
    \subfigure[Effects of $\mu$ and $\sigma_2$ on $\pi_s^*$]{
    \includegraphics[width=7cm]{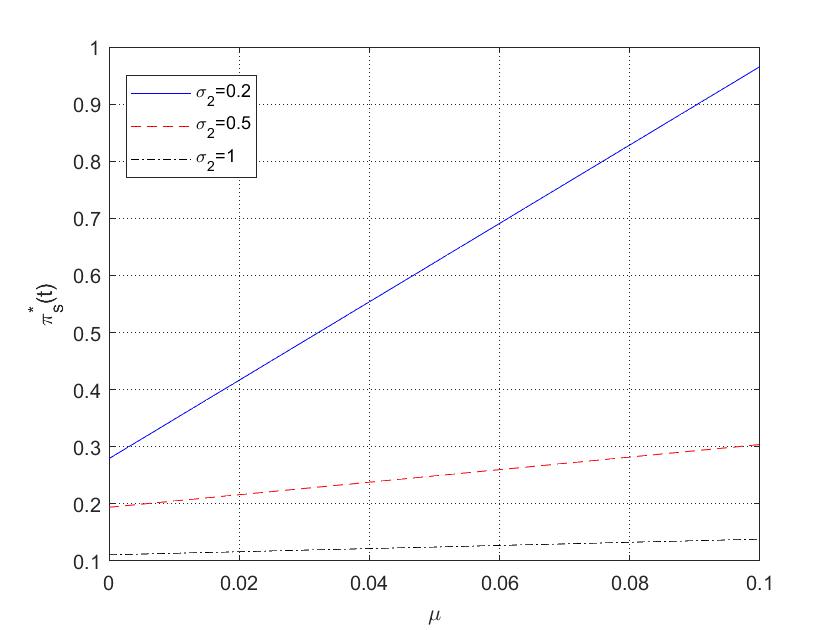}
    }
    \quad
    \subfigure[Effects of $\rho$ and $\sigma_1$ on $\pi_s^*$]{
    \includegraphics[width=7cm]{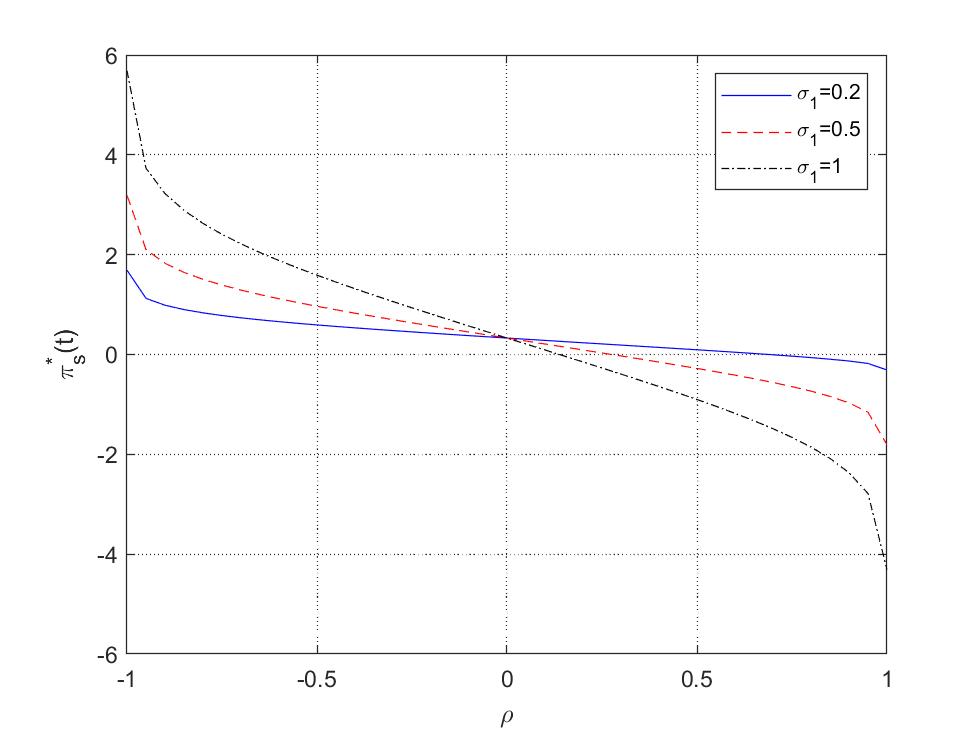}
    }
    \caption{Effects of parameters $\mu$, $\sigma_1$, $\sigma_2$ and $\rho$ on $\pi_s^*$}
    \label{fig:my_label_3}
\end{figure}

\begin{figure}[htpp]
    \centering
    \subfigure[Effects of $t$ and $\mu$ on $\pi_s^*$]{
    \includegraphics[width=7cm]{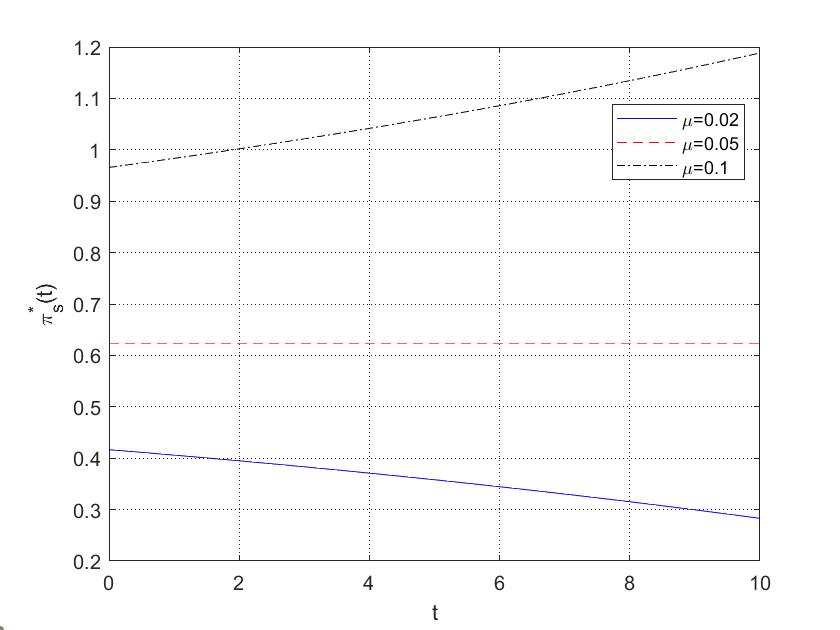}
    }
    \quad
    \subfigure[Effects of $\rho$ and $\gamma$ on $\pi_s^*$]{
    \includegraphics[width=7cm]{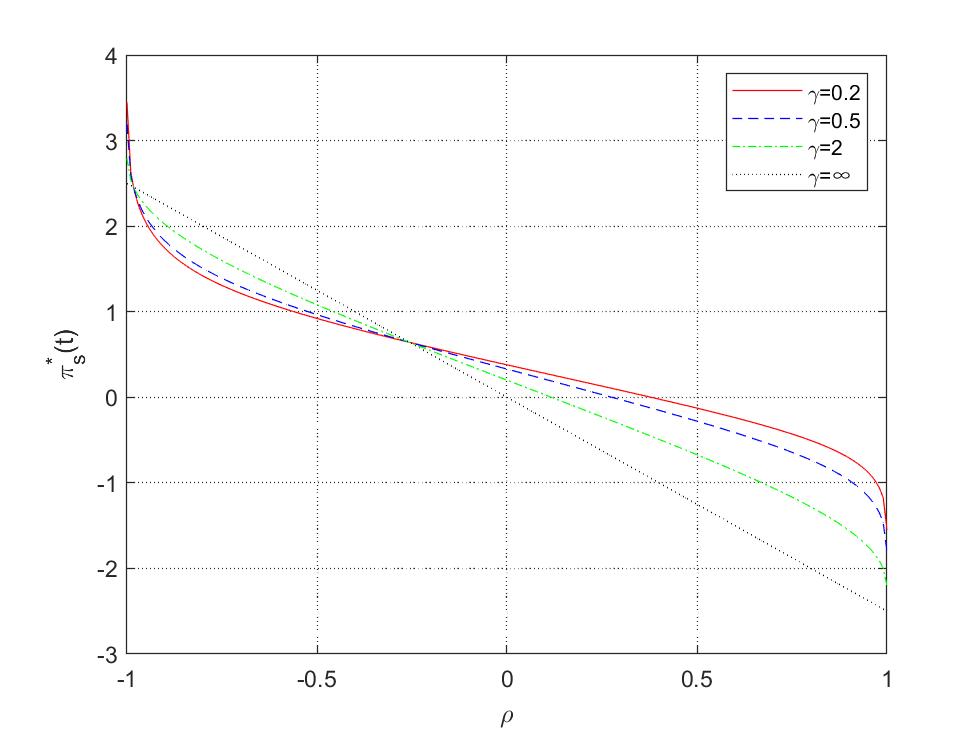}
    }
    \caption{Effects of parameters $t$, $\mu$, $\gamma$ and $\rho$ on $\pi_s^*$}
    \label{fig:my_label_4}
\end{figure}

According to (\ref{436}), we derive that $\frac{\partial \pi_{s}^{*}(t)}{\partial \mu}>0$, $\frac{\partial \pi_{s}^{*}(t)}{\partial t} \geq(\leq$ resp.)0 when $\mu \geq(\leq$ resp.)r, and$\frac{\partial \pi_{s}^{*}(t)}{\partial \sigma_1} \geq(\leq$ resp.)0 when $\rho \leq(\geq$ resp.)0. Figure 3a imply that the optimal investment strategy $\pi_{s}^{*}(t)$ increases with $\mu$, while decreases with $\sigma_2$. This shows that the higher the expected return on stocks, the smaller the volatility, the greater the amount invested in stocks. Figure 3b imply that $\pi_{s}^{*}(t) \geq (\leq$ resp.$) 0$ and is increasing (decreasing resp.) in $\sigma_1$ if $\rho\leq(\geq$ resp.)0. The increase of $\sigma_1$ makes the insurer's surplus more volatile and the insurance business more risky, the insurer will tend to buy more stocks since the stock becomes relatively safer when $\rho \leq 0$. On the contrary, the insurer will tend to short sell more stocks when $\rho \geq 0$. Figure 4a imply that the optimal investment strategy $\pi_{s}^{*}(t)$ increases (decreases resp.) with $t$, when $\mu \geq(\leq$ resp.)r. Obviously, when the expected return on stocks is higher than the risk-free return, the share of investment in stocks increases with $t$. Figure 4b imply that $\pi_{s}^{*}(t)$ increases with $\gamma$ when $-0.9<\rho<-0.4$, decreases with $\gamma$ when $\rho>-0.2$ and $\rho$ very close to -1. The effect of $\rho$ on $\pi_{s}^{*}(t)$ is related to $\beta_1$ and $\beta_2$. The effect of $alpha$, $\beta_1$ and $\beta_2$ have been studied in Bin et al (2016), hence we omit it here.

\begin{figure}[htpp]
    \centering
    \subfigure[Effect of $\delta$ on $\pi_p^*$]{
    \includegraphics[width=7cm]{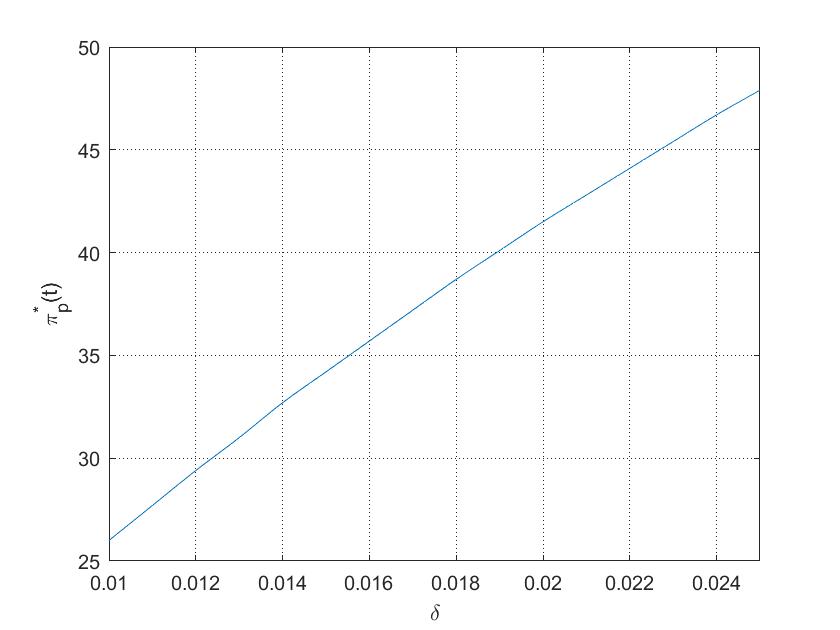}
    }
    \quad
    \subfigure[Effect of $\zeta$ on $\pi_p^*$]{
    \includegraphics[width=7cm]{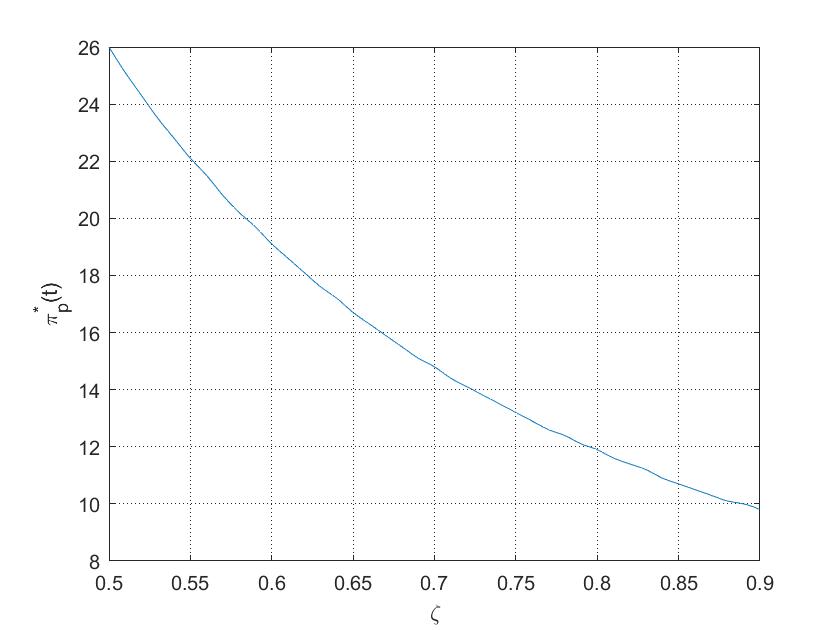}
    }
    \caption{Effects of parameters $\delta$ and $\zeta$ on $\pi_p^*$}
    \label{fig:my_label_5}
\end{figure}

\begin{figure}[htpp]
    \centering
    \subfigure[Effect of $\alpha$ on $\pi_p^*$]{
    \includegraphics[width=7cm]{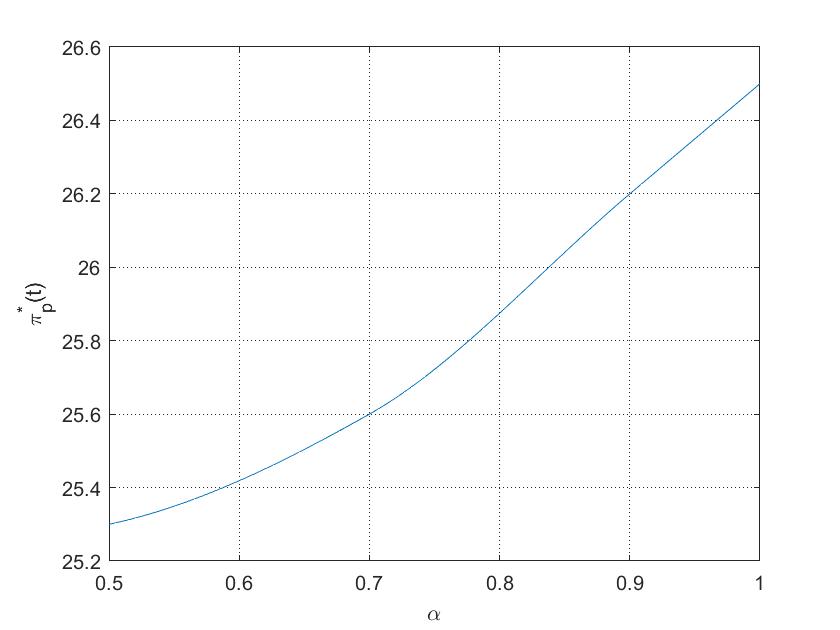}
    }
    \quad
    \subfigure[Effect of $h^{\mathbb{P}}$ on $\pi_p^*$]{
    \includegraphics[width=7cm]{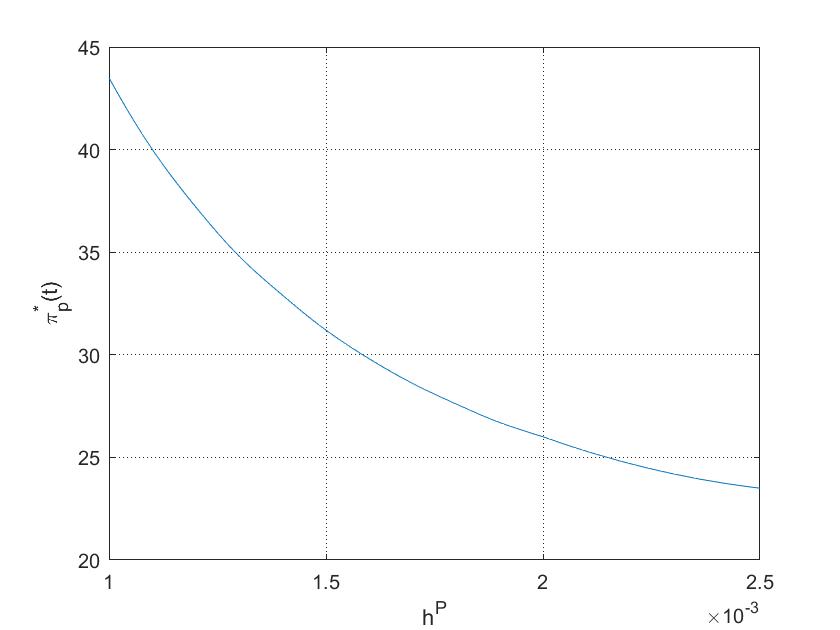}
    }
    \caption{Effects of parameters $\alpha$ and $h^{\mathbb{P}}$ on $\pi_p^*$}
    \label{fig:my_label_6}
\end{figure}

Figure 5 and Figure 6 illustrate the effects of loss rate $\zeta$, credit spread $\delta$, ambiguity attitude $\alpha$ and default intensity $h^{\mathbb{P}}$ on the optimal investment strategy $\pi_p^*$. Figure 5a and 5b imply that the optimal investment strategy $\pi_p^*$ increases with credit spread $\delta$, while decreases with loss rate $\zeta$. The insurer will invest more in the defaultable bond with smaller loss rate $\zeta$ and higher credit spread $\delta$. This finding accords with our financial intuition. The defaultable bond will be more attractive when its default recovery rate $1-\zeta$ is larger and the credit spread is higher.From Figure 6a, we find that $\pi_p^*$ increases with ambiguity attitude $\alpha$. Because defaultable bonds are relatively safe, the more ambiguity-averse insurer will tend to invest more defaultable bond. Observing Figure 6b, we find that the optimal wealth invested in the defaultable bond decreases with the default intensity $h^{\mathbb{P}}$.

\section{Conclusion}
In this paper, we consider the robust optimal reinsurance-investment problem of insurance companies under the $\alpha$-maxmin mean-variance criterion in the defaultable market. Insurance companies can buy reinsurance or acquire new business, and invest in risk-free assets, a stock and a defaultable bond. Since the dynamic mean-variance problem is time-inconsistent, we solve this problem from the perspective of game theory. Using the dynamic programming method, the closed-form expressions and corresponding value functions of the optimal time-consistent strategy in the two cases are obtained. Finally, we analyze the property of the optimal strategy and give the influence of the parameters on the optimal time consistent strategy. From the sensitivity analysis and numerical description, we find that the default event has no effect on the optimal reinsurance strategy and the optimal amount of investment in the stock, while the stock price and the parameters of the insurance risk model affect the insurer’s optimal demand for the defaulted bond. Due to the correlation between stock price and insurance risk, the optimal reinsurance strategy and optimal investment strategy depend on the stock price and the parameters of the insurance risk model. 

In future research, more complex models can be considered, such as the mean variance criterion with state-related risk aversion and the stochastic volatility model with jumps, although it may be difficult to obtain a closed-form solution by doing so. Therefore, other methods can be introduced, such as asymptotic method or other practical methods to deal with the robust optimal reinsurance investment problem.



\newpage
\bibliographystyle{elsarticle-harv}


\end{document}